\documentclass[a4paper,11pt,reqno]{amsart}
\usepackage{amsmath,amssymb,amsthm} 
\usepackage{graphics} 
\usepackage{neuralnetwork}
\usetikzlibrary{shapes}

\usepackage[colorlinks = true]{hyperref}

\usepackage{cleveref}

\usepackage[hyperref,doi=false,url=false,isbn=false, giveninits=true, backref=true,style=alphabetic,maxbibnames=99, backend=biber]{biblatex}
 \bibliography{deep_learning.bib}

 \usepackage{bm}

\numberwithin{equation}{section}

\newtheorem{theorem}{Theorem}[section]
\newtheorem{lemma}[theorem]{Lemma}
\newtheorem{proposition}[theorem]{Proposition}

\newtheorem{corollary}[theorem]{Corollary}

\theoremstyle{definition}

\newtheorem{remark}[theorem]{Remark}

\newcommand{\R}{\mathbb{R}}

\newcommand{\eps}{\varepsilon}

\newcommand{\be}{\begin{equation}}
\newcommand{\ee}{\end{equation}}

\newcommand{\sym}{sym}
\makeatletter
\renewcommand{\fnum@figure}{Fig. \thefigure}
\makeatother
\newcommand{\cost}{\mathsf{c}}

\def\Id{{\rm Id}}
\def\tr{{\rm tr}}

\def\EE{\mathbb{E}}
\def \W{\mathcal{W}}

\newcommand{\cD}{\mathcal{D}}
\newcommand{\cT}{\mathcal{T}}

\newcommand{\bra}[1]{\left( #1 \right)}
\newcommand{\sqa}[1]{\left[ #1 \right]}
\newcommand{\cur}[1]{\left\{ #1 \right\}}
\newcommand{\ang}[1]{\left< #1 \right>}
\newcommand{\abs}[1]{\left| #1 \right|}
\newcommand{\nor}[1]{\left\| #1 \right\|}

\usepackage{enumerate}
\newcommand{\Lip}{\operatorname{Lip}}

\newcommand{\Ker}{\operatorname{Ker}}
\newcommand{\Sym}{\operatorname{Sym}}
\newcommand{\cN}{\mathcal{N}}
\newcommand{\X}{\mathcal{X}}

\newcommand{\cX}{\mathcal{X}}
\newcommand{\cM}{\mathcal{M}}
\newcommand{\cV}{\mathcal{V}}

\newcommand{\Cov}{\operatorname{Cov}}
\title[]{Wide Deep Neural Networks with Gaussian Weights are Very Close to Gaussian Processes }

\author[D. Trevisan]{Dario Trevisan}
\address{D.T.: Dipartimento di Matematica, Università degli Studi di Pisa, 56125 Pisa, Italy. Member of the INdAM GNAMPA group.}
\email{dario.trevisan@unipi.it}

\date{\today  }
\subjclass[2010]{60F05, 60G15, 68T99}
\keywords{Neural Networks, Central Limit Theorem, Wasserstein distance}
\dedicatory{In memory of prof.\ Giuseppe Da Prato}

\begin{document}

\maketitle

\begin{abstract}
We establish novel rates for  the Gaussian approximation of random deep neural networks with Gaussian parameters (weights and biases) and Lipschitz activation functions,  in the wide limit. Our bounds apply for the joint output of a network evaluated any finite input set,  provided a  certain non-degeneracy condition of the infinite-width covariances holds.  We demonstrate that the distance between the network output and the corresponding Gaussian approximation scales inversely with the width of the network, exhibiting faster convergence than the naive heuristic suggested by the central limit theorem. We also apply our bounds to obtain theoretical approximations for the exact Bayesian posterior distribution of the network, when the likelihood is a bounded Lipschitz function of  the network output evaluated on a (finite) training set. This includes popular cases such as the Gaussian likelihood, i.e.\ exponential of minus the mean squared error.
\end{abstract}

\section{Introduction}

In recent years, deep neural networks have emerged as powerful tools in various domains, achieving state-of-the-art performance in tasks such as speech and visual recognition \cite{lecun_deep_2015}, feature extraction and sample generations \cite{goodfellow_generative_2014}. This deep learning revolution has sparked a surge of interest in understanding the scope and limitations of these models, leading to new mathematical investigations in fields like probability, statistics, and statistical physics \cite{roberts2022principles}.

\subsection{Random neural networks} One promising route to analyze deep neural networks is to study their scaling limit  when parameters are  randomly sampled. This approach is motivated by several factors. Firstly, from a Bayesian perspective, a prior distribution on the set of networks is desired, which can then be updated using Bayes' rule based on the observed data, such as training sets in supervised learning. However, in practice, large neural networks are typically trained using iterative optimization algorithms \cite{Goodfellow-et-al-2016}, which often rely on careful initialization of the parameters, often done by random sampling, thus providing a second, more practical, motivation for studying random neural networks. Surprisingly, it has been also observed that even training only the parameters in the last layer of a randomly initialized neural network can yield good performance and significantly accelerate the training process. This intriguing phenomenon has spawned a whole machine learning paradigm \cite{cao2018review}.

In view of these strong motivations, the study of random neural networks has a rich history dating back to pioneers in the field. Early work by Neal \cite{neal_priors_1996} established that wide shallow networks, i.e., networks with a large number of parameters and only one hidden layer, converge to Gaussian processes under suitable tuning of random parameters. This result was derived as a consequence of the central limit theorem. More recent studies \cite{matthews_gaussian_2018, lee_deep_2018} have extended such convergence analysis to deeper networks, defining a path of research that aims to explain the remarkable success of overparameterized neural networks, by establishing a connection with the classical theory of Gaussian processes \cite{williams2006gaussian}.

A significant breakthrough came with the realization that the training dynamics, such as gradient descent, can also be traced in the scaling limit \cite{lee_wide_2019}. This connection leads to a differential equation of gradient flow type associated with the loss functional and the Neural Tangent Kernel, that is the quadratic form associated with the gradient of the parameters evaluated on the training set. Another perspective, known as the mean field limit of networks, see e.g.\ \cite{mei_2018, nguyen_rigorous_2021}, explores the tuning of parameters in a law of large numbers regime, resulting in deterministic nonlinear dynamics for network training.

\subsection{Our contribution}

In light of these developments, our paper contributes to the quantitative analysis of random deep neural networks with Gaussian weights in the wide limit. Our main result, given below as \Cref{thm:main-introduction} in a simplified version, establishes novel rates of convergence for deep networks evaluated on finite input sets and with possibly multidimensional outputs. Our findings shed further light on the behavior of the network output and its distance from the corresponding Gaussian approximation, revealing faster convergence rates than those already established in the literature and (at least to us) somehow unexpected, in view of a naive heuristic suggested by the central limit theorem.

Moreover, we leverage such approximation bounds to obtain analogous rates for the exact Bayesian posterior distribution of the network, when the likelihood is a bounded and Lipschitz function with respect to the network output on a finite training set. This includes the popular choice of Gaussian likelihood, that is the Bayesian counterpart of the mean squared error (see \Cref{sec:bayes} for more details). Thus, we show that the Bayesian point of view on deep neural network is reliable and quantitatively consistent in the wide limit with that provided by Gaussian processes \cite{williams2006gaussian}.

\subsection{Related literature}

To measure our convergence rates, we employ as a metric between probability distributions the Wasserstein distance $\W_p$ (of any given order $p\ge 1$) which is a natural extension of the so-called earth mover's distance in computer science.  Dating back to Monge's works on the optimal transportation of mass, the Wasserstein metric  provides a natural interpretation of the discrepancy between the laws of two random variables, by measuring the most efficient way to give a joint definition in a suitable probability space. Over the past few decades, the theory of optimal transport has experienced significant growth and found applications in various fields beyond calculus of variation and probability, ranging from geometry and partial differential equations \cite{villani2009optimal, ambrosio2005gradient} to statistics and machine learning \cite{peyre2019computational}. In \Cref{sec:wasserstein} we collect all the relevant definitions and basic facts that are needed in the statement of our results. We only anticipate here that to simplify the notation we write $\W_p(X,Y)$ when $X$ and $Y$ are random variables to actually mean the Wasserstein distance of order $p$ between their laws.

The use of a Wasserstein-like distance to quantify the Gaussian approximation of deep neural networks with Gaussian weights can be traced back to \cite{eldan2021non, klukowski2022rate} in the shallow case (i.e., with a single hidden layer), although notably already in the functional framework (i.e., for infinitely many inputs). In the joint work with A.~Basteri, \cite{basteri2022quantitative}, we quantified for the first time the Wasserstein distance of order $2$ between the output laws of a deep neural network with Gaussian weights and a suitable Gaussian process, when evaluated at any finite input set, and established that the distance scales at most as the inverse square root of the hidden layer widths. However, a faster convergence rate was actually observed in numerical simulations \cite[section 6]{basteri2022quantitative}, hinting at the possibility that deep neural networks with random Gaussian weights may be much closer to Gaussian processes. An improvement on the convergence rate was recently established in the case of a single input and one-dimensional output in \cite[Theorem 3.3]{favaro2023quantitative} and thus it was natural to conjecture that an analogue result should hold in the multi-input and multidimensional output case -- this is indeed stated explicitly in \cite[(g) at page 15]{basteri2022quantitative}, however with respect to the convex distance, not the Wasserstein one. Other recent investigations \cite{balasubramanian2023gaussian, cammarota2023quantitative, apollonio2023normal} focus on different distances or on the approximation in functional spaces, rather than improving rates of convergence  for a finite input set. Finally, we mention that exact Bayesian posterior approximation for wide deep neural networks has been qualitatively considered in the literature (see e.g.\ \cite{hron_exact_2020}) but to the author's knowledge this is the first quantitative result in this direction.

\subsection{A simplified statement of our main results}

\begin{figure}
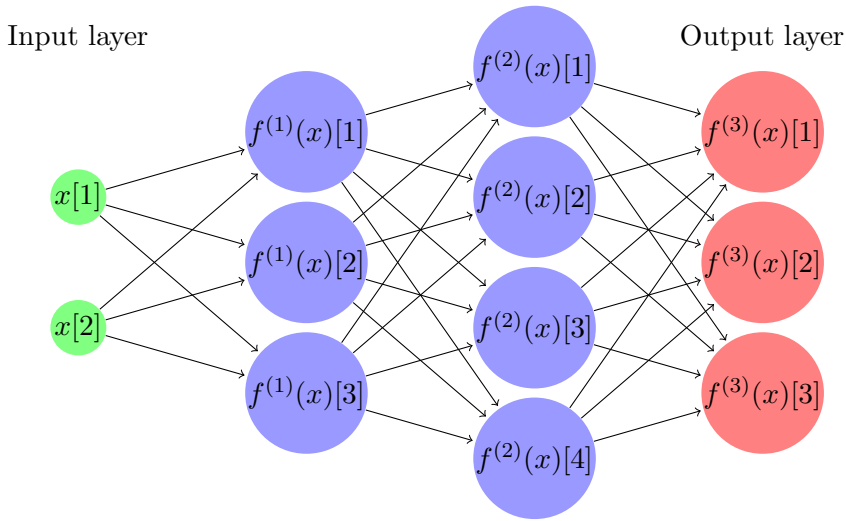

	\begin{neuralnetwork}[height=4, nodesize=17pt, nodespacing=4.5em, layerspacing=3cm, layertitleheight=1em]
		\newcommand{\nodetextclear}[2]{}
		\newcommand{\nodetextx}[2]{$x[#2]$}
		\newcommand{\nodetexthidden}[2]{$f^{(1)}(x)[#2]$}
				\newcommand{\nodetexthiddentwo}[2]{$f^{(2)}(x)[#2]$}
		\newcommand{\nodetexty}[2]{$f^{(3)}(x)[#2]$}
		\inputlayer[count=2, bias=false, title=Input layer, text=\nodetextx]
		\hiddenlayer[count=3, bias=false,  text=\nodetexthidden] \linklayers
		\hiddenlayer[count=4, bias=false,  text=\nodetexthiddentwo] \linklayers
		\outputlayer[count=3, title=Output layer, text=\nodetexty] \linklayers
	\end{neuralnetwork}
	\caption{Graphical representation \cite{cowan_battlesnakeneural_2022} of a fully connected neural network with $L=3$ layers, input size $n_0=2$, output size $n_3=3$ and hidden layer sizes $n_1=3$, $n_2=4$.}\label{fig:nn}
	\end{figure}

In order to make our findings accessible also to an audience that may not be expert in probability theory nor in statistics we give a simplified presentation of our results in the (already interesting) case of a fully connected neural network, where all units (neurons) in each layer are linked with those on the next layer, from the input layer to the output (usually depicted as in \Cref{fig:nn}).

We denote with $L+1$ the total number of layers (including input and output), with sizes $n_0$ for the input layer, $n_L$ for the output layer and $\bm{n} := (n_\ell)_{\ell = 1}^L$ for the hidden and output layers. Then, for every choice of the network parameters, often called weights (matrices) $\bm{W} =(W^{(\ell)})_{\ell=1}^{L}$ and biases (vectors) $\bm{b} = (b^{(\ell)})_{\ell=1}^{L}$, and a (usually fixed) non-linear activation function $\sigma$ acting componentwise on each intermediate layer, a neural network $f^{(L)}: \R^{n_0} \to \R^{n_L}$ is defined.  We refer to \Cref{sec:nn} for a precise definition -- that we actually extend to cover more general architectures.

  For $p \ge 1$, we introduce the constant \begin{equation}\label{eq:gamma-p}  \gamma_p:= ( \sqrt{p}  - \sqrt{p-1})^{1/p} \in [1, \infty) \end{equation} and notice in particular that $\gamma_1 = 1$. Using these concepts and notation, our main result can be stated in a simplified form as follows (we refer to \Cref{thm:main-induction} for a more precise and general version).

 \begin{theorem}\label{thm:main-introduction}
Consider a fully connected neural network $f^{(L)}: \R^{n_0} \to \R^{n_L}$ with Lipschitz activation function $\sigma$ and random
 weights ${\bf W}$ and biases ${\bf b}$ that are independent Gaussian random variables, centered with
 \begin{equation}\label{eq:variances-clean} \EE\sqa{ (W^{(\ell)}_{i,j})^2} = \frac 1 {n_{\ell-1}}, \quad \EE\sqa{ (b^{(\ell)}_{i})^2} = 1, \quad \text{for every $\ell = 1, \ldots, L$ and all $i$, $j$.}\end{equation}
 Then, for every finite set of inputs $\X = \cur{x_i}_{i=1}^k \subseteq \R^{n_0}$ and every $\ell = 1,\ldots, L$, the law of the random variable $f^{(\ell)}[\X] = (f^{(\ell)}(x_i))_{i=1}^k$ is quantitatively close to a vector of $n_{\ell}$ independent and identically distributed Gaussian variables $G^{(\ell)}_{n_{\ell}}[\cX]$, centered and with a (computable) covariance matrix $K^{(\ell)}[\X]$. Precisely, for every $p\ge 1$,
 \begin{itemize}
  \item[a)] if all the covariance matrices $(K^{(\ell)}[\X])_{\ell =1, \ldots, L}$ are invertible, then it holds
  \begin{equation}\label{eq:wp-main-nice}
  \W_p\bra{ f^{(\ell)}[\X],   G^{(\ell)}_{n_{\ell}} [\X]} \le \cost n_\ell \gamma_p^{k n_\ell} \sum_{i=1}^{\ell-1} \frac{1}{n_i},
 \end{equation}
 \item[b)] while, in general, it holds
   \begin{equation}\label{eq:main-worst-case}
  \W_p\bra{ f^{(\ell)}[\X],   G^{(\ell)}_{n_{\ell}} [\X]} \le \cost  \sqrt{n_\ell} \sum_{i=1}^{\ell-1} \frac{1}{\sqrt{n_i}}.
 \end{equation}
 \end{itemize}
In both inequalities the constants\footnote{we use throughout the generic letter $\cost$ to denote a finite positive constant that may depend on some parameters, and may change from line to line.} $\cost< \infty$ do not depend on the layers widths $\bm{n} = (n_i)_{i=1}^{L}$.
 \end{theorem}

The random neural network that we consider is often referred to as the Xavier initialization, where the Gaussian weights are set in \eqref{eq:variances-clean} using a scaling factor that takes into account the number of input neurons in each layer. It is important to note that there are alternative random distributions that can be explored for specific activation functions or network architectures, see e.g.\ \cite{narkhede2022review} for a review.

The \emph{infinite-width} covariance matrices $\bra{ K^{(\ell)}[\X]}_{\ell=1, \ldots, L}$ are uniquely determined by the activation function $\sigma$, the input $\X$ and can be computed recursively, as explained in \Cref{sec:nngp}. In certain cases, such as with ReLU activation, closed formulas can be derived, as discussed in \cite{cho2009kernel, lee_deep_2018}. One notable characteristic is that all $n_\ell$ output components in the Gaussian approximation  $G^{(\ell)}_{n_\ell}[\X]$ are independent and identically distributed variables. While the desirability of this property is often debated, it is important to emphasize that we are considering non-trained networks in this context.

The condition on the invertibility of all the covariance matrices $(K^{(\ell)}[\X])_{\ell =1, \ldots, L}$ is a special case of \emph{non-degeneracy} of the network at the input set introduced in \cite{favaro2023quantitative} (where they discuss higher order analogues as well). On simple examples, such as the ReLU activation, it seems to be generically satisfied. This may explain why the simulations in \cite{basteri2022quantitative} systematically showed better convergence than the ``worst case''  rates, i.e., \eqref{eq:main-worst-case}.

As a consequence of our result and a general measure-theoretic argument based on the dual formulation for the Wasserstein distance of order $1$, we are able to provide bounds for the approximation between the exact Bayesian posteriors associated to a deep random neural network and the corresponding Gaussian process, whenever the likelihood is a  bounded Lipschitz function of the output on a (finite) training dataset. We refer to \Cref{sec:bayes} and in particular \Cref{cor:posterior} for a general presentation of our result. For simplicity, we specialize here the statement to the case of Gaussian likelihood.

 \begin{corollary}
With the same notation of \Cref{thm:main-introduction} for a fully connected neural network $f^{(L)}: \R^{n_0} \to \R^{n_L}$, the associated Gaussian approximation $G^{(L)}_{n_L}: \R^{n_0} \to \R^{n_L}$ and the infinite-width covariances $(K^{(\ell)})_{\ell = 1, \ldots, L}$, consider   any ``training'' dataset $\cur{(x_i, y_i)}_{i\in \cD}$ and a ``test'' set $\cur{x_j}_{j \in \cT}$ (both finite) and define the joint input set $\cX := \cur{x_i}_{i \in \cD} \cup \cur{y_j}_{j \in \cT}$.

Given the Gaussian likelihood function
\begin{equation}
 \mathcal{L} ((z_i)_{i \in \cD} ) := \exp\bra{ - \sum_{i \in \cD} \| z_i - y_i \|^2 }
\end{equation}
write $f^{(L)}[\cX]_{|\cD}$ and respectively $G^{(L)}_{n_L}[\cX]_{|\cD}$ for the Bayesian posterior distributions of the network and respectively the Gaussian process, obtained via Bayes' rule and the likelihood $\mathcal{L}$ evaluated at $z_i:= f^{(L)}(x_i)$ and respectively $z_i:= G^{(L)}_{n_L}(x_i)$, for $i \in \cD$.

Then,
 \begin{itemize}
  \item[a)] if all the covariance matrices $(K^{(\ell)}[\X])_{\ell =1, \ldots, L}$ are invertible, it holds
  \begin{equation}
  \W_1\bra{ f^{(L)}[\X]_{|\cD},   G^{(L)}_{n_{L}} [\X]_{|\cD}} \le \cost \sum_{\ell=1}^{L-1} \frac{1}{n_\ell},
 \end{equation}
 \item[b)] while, in general, it holds
   \begin{equation}
  \W_1\bra{ f^{(L)}[\X]_{|\cD},   G^{(L)}_{n_{L}} [\X]_{|\cD}} \le \cost  \sum_{\ell=1}^{L-1} \frac{1}{\sqrt{n_\ell}}.
 \end{equation}
 \end{itemize}
In both cases the constants $\cost< \infty$ do not depend on the hidden layers widths $\bm{n} = (n_i)_{i=1}^{L-1}$.
 \end{corollary}

As it is well known (see e.g.\ \cite{williams2006gaussian}) for a Gaussian process prior and a Gaussian likelihood, also the posterior law, i.e., $G^{(L)}_{n_{L}} [\X]_{|\cD}$, is Gaussian. Thus, our corollary can be summarized by stating that for a deep neural network with Gaussian prior and Gaussian likelihood, its exact Bayesian predictions are quantitatively close to those of a Gaussian process -- and actually very close if the network is non-degenerate on the entire (test and training) input set.

\subsection{Comments on the proof technique}

The proof of our main result builds upon the fundamental idea, essentially due to \cite{neal_priors_1996}, that the emergence of the Gaussian limit in each layer results from a combination of two factors: firstly, because of the central limit theorem scaling for the weight parameters within each layer, and secondly, because the convergence inherited from the previous layer leads to an almost independence among the neurons. This argument was put forward for deeper architectures in the works \cite{matthews_gaussian_2018, bracale_large-width_2020, hanin_random_2021} using a variety of relatively advanced probabilistic tools, and  was recently made quantitative in \cite{basteri2022quantitative} via a much simpler  argument, arguing via induction and using only elementary properties of the Wasserstein distance. As already mentioned, in \cite[Theorem 3.8]{favaro2023quantitative}, a stronger rate of convergence, comparable to \eqref{eq:wp-main-nice} was first established for a single-input and one-dimensional output of a deep network, under a non-degeneracy assumption of the infinite-width covariances (that are scalars in this case).  Their argument uses Stein's method for Gaussian approximation and appears to be limited to the scalar case, or at least it seems far from obvious how to extend it to the multidimensional case.

Instead, the approach that we put forward in this work still relies on an induction argument over the layers, as in \cite{basteri2022quantitative}, but taking a rather different route. The main idea comes from the problem of studying the properties of an ``empirical kernel'' associated to a neural network $f^{(L)}$ and a smooth function $h: \R^{n_L} \to \R$, defined as the quantity
\begin{equation}
\R^{n_0} \times \R^{n_0} \ni (x,\tilde{x})  \mapsto \frac 1 {n_L} \sum_{i=1}^{n_L} h(f^{(L)}(x)) h(f^{(L)}(\tilde x)).
\end{equation}
Based on the Gaussian approximation and the law of large numbers, one expects that, as the layers widths $(n_i)_{i=1}^L$ becomes large, the following convergence holds:
\begin{equation}\label{eq:convergence-empirical-kernel}
 \frac 1 {n_L} \sum_{i=1}^{n_L} h(f^{(L)}(x)) h(f^{(L)}(\tilde x)) \to \EE\sqa{ h(G^{(L)}(x) ) h(G^{(L))}(\tilde x))},
\end{equation}
where $\bra{ G^{(L)}(x), G^{(L)}(\tilde x)}$ is a Gaussian random variable (with values in $\R^2$) and covariance matrix given by  the infinite-width covariance kernel $K^{(L)}$ evaluated on $\X:=\cur{x,\tilde x}$. Indeed, a simple application of the bounds from \cite{basteri2022quantitative} shows that convergence holds and also provides a quantitative bound (see \eqref{eq:induction-kernel-basteri} below), which however does not reveal Gaussian fluctuations around the limit, as a central limit heuristics would suggest. Our strategy consists precisely in establishing first a precise quantitative description for \eqref{eq:convergence-empirical-kernel}, including the Gaussian fluctuations, and then  somehow  invert the line of reasoning, i.e., starting from the convergence rates for empirical kernel, we obtain our stronger Gaussian approximation rates of the deep neural network. This sketchy description of course does not show why we need the non-degeneracy assumption on the infinite-width covariances, nor the role of several technical results, that are collected in \Cref{sec:technical}.

\subsection{Remarks and open questions}

We highlight further noteworthy aspects of our results and suggest potential directions for future research and extensions in these areas.

\begin{enumerate}
 \item   The dependence on the various parameters of the constants ``$\cost$''  in our bounds is left implicit, to keep the exposition simple, although it seems possible track it and maybe valuable in determining the computational and statistical efficiency of the network under different settings. Specifically, it could be useful to investigate how the constants depend on the depth of the network, the number of neurons per layer, the activation functions used and the dimensionality of the input space.

 \item   It seems feasible, although more technically demanding, to derive stronger rates of convergence also for the derivatives of the network with respect to the input set, by extending the techniques used in our proof. Similarly, for the ``empirical kernels'' built from the derivatives with respect to the parameters (such as the Neural Tangent Kernel), it may be interesting to derive rates of convergence by adapting our analysis. This could also provide a better understanding of the learning dynamics.

 \item   The non-degeneracy condition, while not strictly necessary for convergence, plays a crucial role in obtaining stronger rates of convergence.  It could be relevant to further investigate whether it is indeed necessary for such rates and more importantly establish its genericity in terms of simple properties of the network components.

\item    Our proposed generalized architecture for deep networks (\Cref{sec:nn}) may cover a wider range of network architectures, beyond the fully connected case, e.g.\ including convolutional neural networks and graph neural networks. If this does not follows straightforwardly from our general definition, we believe that the tools from the proof could be then adapted to accommodate different network structures and configurations commonly used in practice.
\end{enumerate}

\subsection{Structure of the paper} \Cref{sec:notation} is devoted to introducing the notation and basic facts about neural networks, Gaussian processes and the Wasserstein distance. \Cref{sec:technical} contains the key technical tools derived from the theory of optimal transportation and quantitative central limit theorems which are utilized in the proof of our main results -- the readers who prefer a less in-depth understanding can safely skip it at first reading. \Cref{sec:proof}  provides a more precise statement of our main result, \Cref{thm:main-induction}, along with its proof. To provide a comparison with the possibly degenerate case, in \Cref{thm:main-induction-basteri} we also extend the bounds from \cite{basteri2022quantitative} to the case of the Wasserstein distance of any order $p \ge 1$.  Lastly, in \Cref{sec:bayes}, we offer a concise overview of the Bayesian approach to supervised learning and outline the application of our results in deriving quantitative bounds for the exact Bayesian posterior distributions.

\subsection{Acknowledgements} The author has been supported by the HPC Italian National Centre for HPC, Big Data and Quantum Computing - Proposal code CN1 CN00000013, CUP I53C22000690001, by the INdAM-GNAMPA project 2023 ``Teoremi Limite per Dinamiche di Discesa Gradiente Stocastica: Convergenza e Generalizzazione'' and by the PRIN 2022 Italian grant 2022WHZ5XH - ``understanding the LEarning process of QUantum Neural networks (LeQun)'', CUP J53D23003890006.

We dedicate this work to the memory of prof.\ Giuseppe Da Prato (1936-2023), distinguished mathematician, who made significant contributions to the study of Gaussian analysis and its applications to evolution equations. Throughout his career,  with his  passion and dedication to rigorous  research he inspired generations of students to explore the field of stochastic analysis.

\section{Notation and basic facts}\label{sec:notation}
We use throughout the letters $S$, $T$ to denote finite sets and $|S|$, $|T|$ to denote the number of their elements. We use the letter $\cost$ to denote some positive and finite constant that may depend on several parameters, e.g.,  $\alpha, \beta, \ldots$, which we can make explicit by writing $\cost = \cost(\alpha, \beta, \ldots)$. We also allow the value of $\cost$ to change from line to line.

\subsection{Tensors} \subsubsection{Vectors and matrices}  We write $\R^{S}$ for the set of  real valued functions $f:S\to \R$, and write  $s \mapsto f(s)$, but also frequently use the notation $f[s] := f(s) = f_s$, to avoid subscripts (e.g.\ in vectors and matrices).
When $S = \cur{1, 2, \ldots, n}$, then we simply write $\R^S = \R^n$, which is the set of $n$-dimensional vectors $v = (v[i])_{i=1}^n$, while if $S = \cur{1, \ldots, m} \times \cur{1, \ldots, n}$ then we write $\R^S = \R^{m\times n}$, the set of matrices $A = (A[i,j])_{i=1, \ldots, m}^{j=1, \ldots, n}$. The space $\R^S$ enjoys a natural structure of vector space with pointwise sum $(f+g)[s] = f[s]+g[s]$ and product $(\lambda f)[s] = \lambda f[s]$ and a canonical basis $(e_s)_{s \in S}$ with $e_s[t] =  1$ if $s=t$ and $e_s[t]=0$ otherwise.

\subsubsection{Tensor products}
Given $A \in \R^{S\times T}$, we use the notation $A[:, t]  := (A[s,t])_{s \in S} \in \R^S$ and similarly $A[s,:] := (A[s,t])_{t \in T} \in \R^T$ (and naturally extended if $A$ belongs to $\R^U$, where $U$ is a product of more than two sets). The tensor product operation between $f \in \R^S$ and $g \in \R^T$ defines $f \otimes g \in \R^{S\times T}$ with $(f \otimes g)[s,t] := f(s) g(t)$. When $f=g$, we simply write $f^{\otimes 2} = f \otimes f$. Any $A \in \R^{S \times T}$ can be identified also with the space of linear transformations from $\R^{T}$ to $\R^S$ via the usual matrix product
\begin{equation} A: f \mapsto Af \quad \text{where} \quad (Af)[s] := \sum_{t \in T} A[s,t] f[t].\end{equation}
The identity matrix $\Id_S \in \R^{S\times S}$ is defined as $\Id_S := \sum_{s} e_s \otimes e_s$, i.e., $\Id_S[s,t] = e_s[t]$, and write $\Id_{n}$ or $\Id_{m\times n}$ in the case of $\R^m$ or $\R^{n\times m}$ respectively.
If $S = S_1 \times S_2$ and $T= T_1 \times T_2$ and $A \in \R^{S_1\times T_1}$, $B \in \R^{S_2 \times T_2}$, then, up to viewing $A \otimes B \in \R^{(S_1\times T_1)\times (S_2 \times T_2)}$ as an element in $\R^{(S_1\times S_2) \times (T_1 \times T_2)}$, the following identity holds:
\begin{equation}
(A \otimes B) (f \otimes g) = (A f) \otimes (B g)
\end{equation}
for every $f \in \R^{T_1}$, $g \in \R^{T_2}$.

\subsubsection{Trace and norms}
We define the trace of $A \in \R^{S\times S}$ as
\begin{equation}
 \operatorname{tr}(A) := \sum_{s\in S} A[s,s].
\end{equation}
 We  use this notion to define a scalar product on $\R^S$ as
 \begin{equation}
 \ang{f, g} := \operatorname{tr}(f \otimes g) = \sum_{s\in S} f[s]g[s],
 \end{equation}
which induces the norm
\begin{equation}\nor{f} =\sqrt{ \operatorname{tr}( f^{\otimes 2})} = \sqrt{ \sum_{s \in S} |f[s]|^2},\end{equation} yielding the Euclidean norm in case of a vectors in $\R^k$ and the Frobenius (or Hilbert-Schmidt) norm for  matrices $\R^{n\times m}$. The tensor product is associative and the norm is multiplicative, i.e., $\nor{f\otimes g} = \nor{f}\nor{g}$.

The induced operator norm of a linear operator $A \in \R^{T \times S}$ will be denoted instead as
\begin{equation}
 \nor{A}_{op} := \sup_{ \substack{ v \in \R^S \\ \nor{v} = 1}} \nor{ Av},
\end{equation}
It holds
\begin{equation}\label{eq:norm-hs-op}
 \nor{A}_{op} \le \nor{A} \le \sqrt{|S|} \nor{A}_{op}.
\end{equation}

\subsubsection{Lipschitz functions and polynomial growth}

We say that a function $f: \R^S \to \R^T$ has polynomial growth of order $q >0$ if, for some constant that we denote with $\nor{f}_q< \infty$ it holds
\begin{equation}
  \nor{ f(x) } \le \nor{f}_q \bra{1+ \nor{x}^q} \quad \text{for every $x \in \R^S$.}
 \end{equation}
A function $f: \R^S \to \R^T$ is said to be Lipschitz continuous if there exists a constant $\Lip(h) < \infty$ such that $\nor{h(x) -  h(y)} \le \Lip(h) \nor{x-y}$ for every $x, y \in \R^S$. Notice that a Lipschitz continuous function has polynomial growth of order $q=1$, since for every $x \in \R^S$, we have
\begin{equation}\label{eq:linear-growth}\begin{split}
 \nor{h(x)} & \le \nor{h(0)} + \nor{h(x) - h(0)} \le \nor{h(0)} + \Lip(h) \nor{x} \\
 & \le \nor{h}_{\Lip} (1+\nor{x})
\end{split}\end{equation}
where we set $\nor{h}_{\Lip} :=  \nor{h(0)}+ \Lip(h)$.

\subsubsection{Symmetric operators}
Any $A \in \R^{S \times S}$ can be seen as a quadratic form over $\R^S$: it is  symmetric if $A$ equals its transpose, $A^\tau[s,t] := A[t,s]$, i.e., $A[s,t]:=A[t,s]$ for $s$, $t \in S$. A symmetric $A \in \R^{S\times S}$ is positive semidefinite if for every $v \in \R^S$, it holds
\begin{equation}
  \sum_{s,t \in S} v[s] A[s,t] v[t]  \ge 0.
\end{equation}
We write $\Sym^S$ (respectively, $\Sym_+^S$) for the set of symmetric  (respectively, symmetric and positive semidefinite) $A \in \R^{S\times S}$ and simply $\Sym^k$ (respectively, $\Sym_+^k$) if $S = \cur{1, \ldots, k}$.
 Given $A \in \Sym_+^S$ and $B \in \Sym^S$, the following notation will be useful: we write
 \begin{equation}\label{eq:ac-matrices}
  B \ll A \quad \text{ if and only if } \quad \text{for every $v \in \R^S$ such that $Av =0$, it holds $Bv = 0$,}
 \end{equation}
 which can be seen as a non-commutative analogue of absolute continuity between measures. Notice that $B \ll A$ trivially holds for every $B$ if $A \in \Sym_+^S$ is invertible.

If $A \in \Sym^S_+$, we write $\lambda(A)$ for the smallest  eigenvalue of $A$ (which may be null) and $\lambda_+(A)$ for the smallest strictly positive eigenvalue of $A$, with the convention that $\lambda_+(0) = \infty$ if $A = 0$, so expressions like $\nor{A} / \lambda_+(A)$ are also well-defined (we conventionally let $0/\infty = 0$).

\subsubsection{Square roots}
Given  $A\in \Sym_+^S$, its (positive) square root $\sqrt{A}$ is defined as usual via spectral theorem as $\sqrt{A} = U \sqrt{D} U^\tau$ where $U \in \R^{S\times S}$ is orthogonal, i.e., $U U^\tau = \Id$ and $U^\tau A U = D= \sum_{i=1}^m \lambda_i e_i \otimes e_i$ is diagonal with $\lambda_i \ge 0$ (the eigenvalues of $A$) and $\sqrt{D} = \sum_{i=1}^m \sqrt{\lambda_i}e_i \otimes e_i$ is the diagonal with entries being the square roots of the entries of $D$ (and $m = |S|$). For $A \in \R^{S \times S}$ (not  necessarily symmetric), we write
\begin{equation}\label{eq:abs-A}
 |A| := \sqrt{ A A^\tau},
\end{equation}
where $A^\tau$ is the transpose of $A$.

The map $A \mapsto \sqrt{A}$ is H\"older continuous with exponent $1/2$ with respect to the operator norms:
\begin{equation}\label{eq:holder}
 \nor{ \sqrt{A} - \sqrt{B} }_{op} \le \cost \sqrt{\nor{ A - B}_{op}} \quad \text{for every $A$, $B \in \Sym_+^S$.}
\end{equation}
for some constant $\cost < \infty$. In view of \eqref{eq:norm-hs-op}, we also have
\begin{equation}\label{eq:holder}
 \nor{ \sqrt{A} - \sqrt{B} } \le  \cost \sqrt{\nor{ A - B}} \quad \text{for every $A$, $B \in \Sym_+^S$.}
\end{equation}
with $\cost = \cost(S)< \infty$.

As for the case of scalars, the map  $\sqrt{\cdot}$ is smooth away from zero. Precisely, we use the following facts: given $A \in \Sym_+^S$ and $B \in \Sym^S_+$ such that $B \ll A$ (in the sense of \eqref{eq:ac-matrices}), then \cite[Proposition 3.2]{van1980inequality} gives
\begin{equation}\label{eq:ando}
\nor{ \sqrt{A} - \sqrt{B}} \le \frac{1}{\sqrt{ \lambda_+(A) } } \nor{ A - B} \le \cost \nor{A-B},
\end{equation}
where $\cost = \cost(S, A)<\infty$. If moreover $\nor{B}_{op}\le  \lambda(A)/2$, then
\begin{equation}\label{eq:taylor-sqrt}
   \nor{ \sqrt{ A+B }  - \sqrt{A} - D_{\sqrt{\cdot}}(A) B} \le \cost \frac{ \nor{ B}^2 }{ \lambda_+(A)^{3/2}} \le \cost \nor{B}^2,
\end{equation}
where again we let  $\cost = \cost(S, A)< \infty$ and the map $B \mapsto D_{\sqrt{\cdot}}(A)B \in \Sym^S$ linear on the subspace of $B \in \Sym_+^S$ such that $B \ll A$ and
\begin{equation}\label{eq:norm-dsqrta}
 \nor{D_{\sqrt{\cdot}}(A) B} \le \cost  \frac{\nor{B}}{\lambda_+(A)^{1/2}} \le \cost \nor{B},
\end{equation}
where $\cost = \cost(S,A)<\infty$.
When $A = \Id_S$, it is not difficult to show that  $D_{\sqrt{\cdot}}(A) B = B/2$. Finally, we remark again that the condition $B \ll A$ is void if $A$ is invertible.

\subsubsection{Partial trace}

Given  $B\in \R^{(S \times T) \times (S \times T)}$, its partial  trace with respect to $S$ is defined as $\tr_{S}(B) \in \R^{T \times T}$, given by the expression
\begin{equation}
 \tr_{S}(B) [t_1, t_2] := \sum_{s \in S} B[s,t_1,s,t_2].
\end{equation}
Clearly, $B \mapsto \tr_S(B)$ is linear with
\begin{equation}\label{eq:norm-partial-trace}
 \nor{ \tr_S(B) } \le \sqrt{|S|} \nor{ B},
\end{equation}
and moreover it is not difficult to check that $\tr_{S}(B) \in \Sym_+^T$ if $B \in \Sym_+^{S\times T}$.  If $B = A^{\otimes 2}$, we then write
\begin{equation}\label{eq:abs-A-s}
 | A |_S =  \sqrt{\tr_{S}(A^{\otimes 2})  }.
\end{equation}

When $S = \cur{1,\ldots, n}$, we simply write $\tr_n(B) := \tr_S(B)$ and $|A|_n = |A|_S$. We warn the reader that the notation $\tr_S$ and $\tr_n$ can be slightly ambiguous if more copies of $S$ (or $n$) are considered -- for example if $S=T$. We try however to be consistent with the notation and highlight below possible inconveniences.  Let us notice for example that if $S = T$, then $|A|_S = |A| = \sqrt{A A^\tau}$ as defined in \eqref{eq:abs-A}.

\subsection{Random variables}

\subsubsection{Laws and moments} Given a random variable $X$ with values in $\R^S$, defined on some probability space  $(\Omega, \mathcal{A}, \mathbb{P})$, we write $\mathbb{P}_X$ for its law, i.e., $\mathbb{P}_X(E) = \mathbb{P}(X \in E)$ for $E \subseteq \R^S$ Borel. We write $X \stackrel{law}{=} Y$ if two random variables $X$, $Y$ have the same law.
For $p\ge 1$, the absolute moment of order $p$ of a random variable $X$ with values on $\R^S$ is defined as the quantity $\EE\sqa{ \nor{X}^p}$.  We also write $\nor{X}_{L^p} := \EE\sqa{ \nor{X}^p }^{1/p}$ for the usual Lebesgue norm of order $p$.

\subsubsection{Mean and covariance} We write $\EE\sqa{X} \in \R^S$ for the  (componentwise) mean value of $X$, i.e., $\EE\sqa{X}[s] = \EE\sqa{X[s]}$, which is well-defined provided that $X$ has finite absolute moment of order $1$, and  $\EE\sqa{ X^{\otimes 2}}$ for its second moment, i.e., the collection $(\EE\sqa{X[s] X[t]})_{s,t\in S}$ (provided that $X$ has finite absolute moment of order $2$). It is simple to prove that $\EE\sqa{ X^{\otimes 2}} \in \Sym_+^S$ and that a.s.\ it holds
\begin{equation}
 \label{eq:Xtwo-ac}
  X^{\otimes 2} \ll \EE\sqa{ X^{\otimes 2}}.
\end{equation}
Indeed, if $v \in \R^S$ is such that $\EE\sqa{ X^{\otimes 2}} v = 0$, then $\EE\sqa{ \abs{\ang{X, v}}^2} = \ang{v,\EE\sqa{ X^{\otimes 2} }v } = 0$, hence a.s.\ $\ang{X, v} = 0$ and therefore $X^{\otimes 2} v   = X \ang{X, v} =0$ as well.

The covariance $\Cov(X) \in \Sym_+^S$ of $X$ is defined as the second moment of the centered variable $X-\EE\sqa{X}$, which coincides with its second moment if $X$ is centered, i.e., $\EE\sqa{X} =0$. Recalling that $\nor{X}^2 = \tr(X^{\otimes 2})$, and exchanging expectation with the trace operation, we have that
\begin{equation}
 \nor{X}_{L^2}^2 = \tr( \EE\sqa{X^{\otimes 2}} ),
\end{equation}
 which further equals $\tr( \Cov(X))$ if $X$ is centered.

\subsubsection{Gaussian variables}
Given $K \in \Sym_+^S$, we use the notation $\cN_K$ to denote a random variable with values on $\R^S$  with centered Gaussian distribution and covariance $K$, so that
\begin{equation}
 \EE\sqa{ \cN_K} = 0 \quad \text{and } \quad \Cov(\cN_K) = \EE\sqa{ \cN_K^{\otimes 2}} = K.
\end{equation}
if $K = \Id_{S}$, we write $\cN_S:=\cN_{\Id_S}$, if $S = k$ we write $\cN_k=\cN_{\Id_k}$, and if $S = n \times k$ we write $\cN_{n \times k} = \cN_{\Id_{n\times k}}$. Although this notation could be slightly ambiguous since we are not specifying the probability space where $\cN_K$ is defined, this will eventually be quite handy because in most cases $\cN_K$ will be independent of (most of) the variables appearing the in expressions. For example, we have the identity between the laws
\begin{equation}
 \cN_K \stackrel{law}{=} \sqrt{K} \cN_S.
\end{equation}
Let us recall that for every $p \ge 1$ there exists a constant $\cost = \cost(p, S) < \infty$ such that, for any Gaussian random variable $\cN_K$ with values in $\R^S$ it holds
\begin{equation}\label{eq:moment-p-gaussian}
 \nor{ \cN_K}_{L^p} \le \cost  \sqrt{\tr\bra{K }} = \cost \nor{ \sqrt{K}}.
\end{equation}

The notation is $\cN_K$ is actually so useful that we naturally extended to the case where $K$ is a random variable itself (taking values in $\Sym_+^S$): in such a case it  means by definition the variable
\begin{equation}
 \cN_K := \sqrt{K} \cN_S,
\end{equation}
where $\cN_S$ is a standard Gaussian variable independent of $K$ (possibly defined by enlarging the original probability space where $K$ is originally defined). Equivalently,  $\cN_K$ can be defined as a random variable whose conditional law with respect to the value $K=\kappa \in \Sym^S_+$ is Gaussian centered with covariance $\kappa$. In particular, if $A$ is any random variable with values in $\R^{T \times S}$, independent of a standard Gaussian variable $\cN_S$, we have the identity in law
\begin{equation}
 A \cN_S \stackrel{law}{=} \cN_{ \tr_S (A^{\otimes 2})  } =   |A|_S \cN_T.
\end{equation}
where   $|A|_{S}  \in \Sym^{S}_+$ is defined according to \eqref{eq:abs-A-s}. Using this fact we collect for later use the following  identity between laws: if $A$ is a random variable taking values in $\R^{T_1 \times S}$ and $\cN_{T_2 \times T_1}$ is independent of $A$, then
\begin{equation}\label{eq:identity-in-law}
  (\cN_{T_2 \times T_1} \otimes \Id_S ) A  \stackrel{law}{=} \mathcal{N}_{\Id_{T_2} \otimes  \tr_{T_1}(A^{\otimes 2}) } \stackrel{law}{=} \bra{\Id_{T_2 } \otimes   |A|_{T_1}} \cN_{T_2 \times S}.
\end{equation}
 Indeed, the three expressions above define random variables whose conditional law with respect to $A$ is (conditionally) Gaussian centered with covariance $\Id_{T_2} \otimes \tr_{T_1}(A^{\otimes 2})$.

 Finally, we also need the following observation: if $\cN_A$ is a Gaussian random variable with values in $\R^{S \times T}$, so that $A \in \Sym_+^{(S \times T) \times (S\times T)}$ then, by linearity of the partial trace, $\tr_S(\cN_A)$ is also Gaussian random variable with values in $\R^T$. We write in what follows:
 \begin{equation}\label{eq:partial-trace-SS}
   \tr_{S\times S} (A) := \Cov\bra{ \tr_S(\cN_A)} \in \Sym_+^{T \times T}.
 \end{equation}
 This is justified by  the following straightforward identity, valid for $t_1, t_2, \tilde{t}_1, \tilde{t}_2 \in T$:
 \begin{equation}
  \Cov\bra{ \tr_S(N_A)} [(t_1, t_2), (\tilde t_1, \tilde t_2)] = \sum_{s, \tilde s \in S} A[ (s, t_1, s, t_2), (\tilde{s}, \tilde{t_1},\tilde{s}, \tilde{t_2} )],
 \end{equation}
We remark however that the notation $\tr_{S\times S}$ could be regarded as ambiguous, as we are not precisely specifying which pairs of indices on $S$ are traced out.

\subsection{Wasserstein distance}\label{sec:wasserstein}

\subsubsection{Notation} We recall, mostly without proof, some known facts about optimal transportation theory, referring to any monograph such as \cite{villani2009optimal, ambrosio2005gradient, peyre2019computational, santambrogio2015optimal} for details. The  Wasserstein distance of order $p\ge 1$ between  two probability measures $\mu$, $\nu$ on $\R^S$  is defined  as
\begin{equation}\label{eq:wp} \W_p(\mu, \nu) := \inf\cur{ \nor{X-Y}_{L^p} \, :  \text{$X$, $Y$ random variables with $\mathbb{P}_X = \mu$, $\mathbb{P}_Y = \nu$}},\end{equation}
where the infimum runs over all the random variables $X$, $Y$ (jointly defined on some probability space) with marginal laws  $\mu = \mathbb{P}_X$ and $\nu=\mathbb{P}_Y$ (also called couplings).
For simplicity, we allow for a slight abuse of notation in all what follows and write
\begin{equation}
 \W_p(X,Y) := \W_p(\mathbb{P}_X, \mathbb{P}_Y),
\end{equation}
whenever we consider two random variables $X$, $Y$ with values on $\R^S$. Notice that it does not matter whether $X$ and $Y$ are already  defined on the same probability space, since $\W_p(X, Y)$ still is defined in terms of their marginal laws (and minimization runs over all the possible joint definitions of $X$ and $Y$).

\subsubsection{General properties}
If however $X$ and $Y$ are defined on the same probability space, then it trivially holds
\begin{equation}\label{eq:w-2-trivial-bound} \W_p(X,Y) \le  \nor{X-Y}_{L^p}.\end{equation}

In order to ensure that the quantity $\W_p(X,Y)$ is finite and it defines an actual distance, we (tacitly) restrict our considerations to random variables with finite absolute moments of order $p$. Then, a  sequence of random variables $(X_n)_n$ converges towards a random variable $X$, i.e., $\lim_{n \to \infty}\W_p(X_n, X) = 0$ if and only if $\lim_{n \to \infty} X_n$ in law together with their absolute moments of order $p$, i.e., $\lim_{n \to \infty} \nor{X_n}_{L^p}  = \nor{X}_{L^p}$.
Given random variables $X$, $Y$, $Z$ with values $\R^S$, the triangle inequality reads
\begin{equation}\label{eq:triangle}
\W_p(X,Z) \le \W_p(X,Y)+\W_p(Y,Z),
\end{equation}
and, if $Z$ is independent of $X$ and $Y$,  the following inequality holds:
\begin{equation}\label{eq:sum} \W_p(X+Z, Y+Z) \le \W_p(X, Y).\end{equation}

The distance $\W_p$ raised to the power $p$ (also known as the transportation cost) enjoys a useful convexity property, which we can state as follows: given independent random variables $X_1$, $X_2$, $U$ with values respectively in $\R^{S_1}$, $\R^{S_2}$, $\R^V$ and Borel functions $\varphi_1:\R^{S_1}\times \R^{V} \to \R^{T}$, $\varphi_2: \R^{S_2} \times \R^{V} \to \R^T$, it holds
\begin{equation}\label{eq:convexity} \W_p\bra{\varphi_1(X_1, U), \varphi_2(X_2, U) } \le  \bra{ \int_{\R^V} \bra{\W_p\bra{\varphi_1(X_1, u), \varphi_2(X_2, u) }}^p  d \mathbb{P}_U(u)}^{1/p}.\end{equation}

\subsubsection{Optimal transport between Gaussian laws}

If $X$, $Y$ are Gaussian random variables with values on $\R^S$, then it holds
\begin{equation}\label{eq:gaussian} \W_p(X, Y) \le \cost \bra{ \nor{\EE\sqa{X} - \EE\sqa{Y}} + \nor{ \sqrt{\Cov(X)} - \sqrt{\Cov(Y)}}},\end{equation}
for some constant $\cost = \cost(p, S)< \infty$. This can be easily seen by considering a standard Gaussian variable $\cN_S$  and letting $X = \EE\sqa{X}+ \sqrt{\Cov(X)} Z$ and $Y = \EE\sqa{Y}+ \sqrt{\Cov(Y)}Z$. Then, by \eqref{eq:w-2-trivial-bound} and \eqref{eq:moment-p-gaussian},
\begin{equation}
\begin{split}
\W_p(X, Y)  & \le  \nor{ \bra{\EE\sqa{X}-\EE\sqa{Y}} +  \bra{ \sqrt{ \Cov(X) } - \sqrt{\Cov(Y)} } \cN_S }_{L^p}\\
& \le  \cost  \bra{ \nor{ \EE\sqa{X}-\EE\sqa{Y} } +  \nor{\sqrt{ \Cov(X) } - \sqrt{\Cov(Y)}}}.
\end{split}
\end{equation}

\subsection{Deep Neural Networks}\label{sec:nn}

\subsubsection{Fully connected architecture}
The feed-forward fully connected network (also known as multi-layer perceptron) is the simplest deep neural network architecture and consists of a sequence of hidden layers stacked between the input and output layers, where each node in a layer  is connected with all nodes in the subsequent layer. Let $L \ge 1$ denote the total number of layers (excluding the input one), let $d_0$ denote the dimension of the input space,
\begin{equation}
 \bm{n} = (n_1, \ldots, n_L)
\end{equation}
denote the size of the layers, so that network output is a vector in $\R^{n_{L}}$. The construction goes as follows: we fix a family of functions (called activation functions)
\begin{equation}
 \bm{\sigma} = (\sigma^{(1)}, \sigma^{(2)},\ldots, \sigma^{(L)} ) \quad \text{with} \quad \sigma^{(\ell)}: \R \to \R,
\end{equation}
and two  sets of parameters (called respectively weights and biases),
\begin{equation}
\bm{W} = (W^{(1)}, W^{(2)}, \ldots, W^{(L)}), \quad  \bm{b} = (b^{(1)}, b^{(1)}, \ldots, b^{(L)}),
\end{equation}
where, for every $\ell =1, \ldots, L$, one has $W^{(\ell)} \in \R^{n_{\ell} \times n_{\ell-1}}$ and $b^{(\ell)} \in \R^{n_{\ell}}$. We begin with defining $n_0:= d_0$ and
\begin{equation}f^{(0)} : \R^{d_0} \to \R^{n_0}, \quad f^{(0)}(x) = x,\end{equation}
and recursively, for $\ell=1, \ldots, L$,
\begin{equation}\label{eq:fully-connected} f^{(\ell)} : \R^{d_0} \to \R^{n_{\ell}}, \quad f^{(\ell)}(x) = W^{(\ell)}\sigma^{(\ell)} (f^{(\ell-1)}(x)) + b^{(\ell)}.
 \end{equation}
The composition with $\sigma^{(\ell)}$ is always to be understood componentwise, i.e.,
\begin{equation}
 \sigma^{(\ell)} (f^{(\ell-1)}(x)) [i] := \sigma^{(\ell)}( f^{(\ell-1)}(x)[i]).
\end{equation}
It is quite common to choose $\sigma^{(1)}(x) = x$, the identity function. Popular choices for $\sigma^{(\ell)}$ include the function $\sigma(x) = \max\cur{x, 0}$ (called ReLU).

\subsubsection{Generalized architecture}
In order to deal efficiently with the outputs of a network evaluated at an input set $\X = \cur{x_i}_{i=1}^k \subseteq \R^{d_0}$, it is convenient to slightly generalize the above construction, allowing for more general activation functions and larger output dimensions. Precisely, we fix a family of auxiliary dimensions
\begin{equation}
 \bm{S} = ( S_1, \ldots,  S_L) \quad \text{and} \quad \bm{a} = (a_1, \ldots, a_L)
\end{equation}
where $S_i$ are (finite) sets, and extend the notion of activation functions $\bm{\sigma}$ allowing for
\begin{equation}\label{eq:sigma-gen}
 \sigma^{(\ell)}: \R^{S_{\ell-1}} \to \R^{a_\ell \times S_{\ell} },
\end{equation}
where it will be useful to set $S_0: = d_0 \times k$, where $k$ denotes the size of the input set $\cX$. Correspondingly, we group together the weights and biases in a single family of parameters, that we still denote with $\bm{W} = (W^{(1)}, \ldots, W^{(L)})$ where now $W^{(\ell)} \in \R^{n_{\ell} \times (n_{\ell-1} \times a_{\ell}) }$ and  we conveniently let $n_0:=1$.

With this notation, we now define $f^{(\ell)} \in \R^{n_\ell \times  S_\ell}$ recursively, by letting
\begin{equation}
 f^{(0)} :=  \sum_{i=1}^k x_i e_i \in \R^{d_0 \times k} = \R^{1\times S_0 },
\end{equation}
i.e., we identify the input set $\X$ with a  matrix
and, for $\ell=1, \ldots, L$, we let
\begin{eqnarray}\label{eq:f-induction}
f^{(\ell)}  := (W^{(\ell)} \otimes \Id_{S_\ell})  \sigma^{(\ell)} ( f^{(\ell-1)} )
 \end{eqnarray}
 where $\sigma^{(\ell)}$ acts componentwise, i.e., writing
 \begin{equation}
  f^{(\ell-1)} = (f^{(\ell-1)}[i,:])_{i=1,\ldots, n_{\ell-1}}
 \end{equation}
 with $f^{(\ell-1)}[i,:] \in \R^{S_{\ell-1}}$, we set
 \begin{equation}
  \sigma^{(\ell)}(f^{(\ell-1)}) := (\sigma^{(\ell)} ( f^{(\ell-1)}[i,:] ))_{i=1,\ldots, n_{\ell-1}}.
 \end{equation}

The  fully connected architecture when evaluated at an input set $\X = \cur{x_i}_{i=1}^k$ can be then recovered by setting  $S_0 := d_0 \times k$, $S_{\ell} := k$ for every $\ell=1,\ldots, L$  and $a_{1} := d_0 +1$ and $a_{\ell} = 2$ for $\ell=2, \ldots, L$. Then, for every $v \in \R^{S_0} = \R^{d_0 \times k}$ we define
\begin{equation}
 \sigma^{(1)}( v ) = (v, 1_k)^\tau \in \R^{(d_0+1) \times k},
\end{equation}
where $1_k \in \R^k$ denotes a vector of length $k$ with $1_k[i] = 1$ for every $i=1, \ldots, k$, and for $\ell=2,\ldots L$, $v\in \R^{k}$, we set
\begin{equation}
 \sigma^{(\ell)}(v) = ( \sigma(v), 1_k)^{\tau} \in \R^{2 \times k},
\end{equation}
where $\sigma(v) = (\sigma(v[i]))_{i=1, \ldots, k}$ is applied componentwise. The additional vector $1_k$ takes into account the bias term and it is a straightforward computation to realize that with these definitions the resulting $f^{(\ell)}$ recovers precisely the multi-layer perceptron architecture evaluated at the input set $\cX$.

\subsubsection{Neural Network Gaussian Processes} \label{sec:nngp}
The Gaussian process approximating a wide random fully connected neural network is defined  \cite{lee_deep_2018} in terms of its covariance operator and essentially depends only on the choice of the activation function $\sigma$.
Instead of describing the specific construction stemming from the multilayer perceptron architecture \eqref{eq:fully-connected},  it is  more convenient to give the definition for the generalized architecture we introduced above. For a given sequence of auxiliary dimensions $\bm{S}$, $\bm{a}$ and activation functions $\bm{\sigma}$ as in \eqref{eq:sigma-gen} (with $S_0 :=d_0 \times k$ the size of the given input set $\X = \cur{x_i}_{i=1}^k$ and $n_0 :=1$), we begin with conveniently defining $G^{(0)} \in \R^{S_0}$ as the random variable that is constantly equal to $f^{(0)}$ i.e.,
\begin{equation}
 G^{(0)} = \sum_{i=1}^k x_i e_i \in \R^{d_0 \times k} = \R^{1 \times S_0},
\end{equation}
so that in particular its covariance is null, $K^{(0)} = 0 \in\Sym_+^{S_0}$. For $\ell=1, \ldots, L$, having already defined a Gaussian variable $G^{(\ell-1)}$ with values in $\R^{k_{\ell-1} }$ we introduce the \emph{infinite-width} covariance
\begin{equation}\label{eq:k-ell}
  K^{(\ell)}  := \tr_{a_\ell}\bra{\EE\sqa{   \bra{\sigma^{(\ell)}( G^{(\ell-1)}}^{\otimes 2} } } \in \Sym_+^{S_\ell}
  \end{equation}
  and define the associated Gaussian random variable
  \begin{equation}
   G^{(\ell)} := \cN_{K^{(\ell)}}.
  \end{equation}
 Notice that this construction yields a sequence of Gaussian random variables $G^{(\ell)}$  with values in $\R^{S_\ell}$, which is  different from the space of values taken by a network $f^{(\ell)}$, i.e. $\R^{n_\ell \times S_\ell}$. Indeed, the random variables actually approximating $f^{(\ell)}$, that we denote with $G^{(\ell)}_{n_\ell}$, are simply defined as a collection of $n_\ell$ independent copies of $G^{(\ell)}$, so that they have Gaussian  law with covariance
\begin{equation}
 \Cov( G^{(\ell)}_{n_\ell}) :=  \Id_{n_\ell} \otimes K^{(\ell)},
\end{equation}
i.e., $G^{(\ell)}_{n_\ell} = \cN_{\Id_{n_\ell} \otimes K^{(\ell)} }$.

Finally, we borrow from \cite{favaro2023quantitative} the following definition: we say that the infinite-width covariances $\bm{K} := (K^{(\ell)})_{\ell=1, \ldots, L}$ are non degenerate (on the input set $\X$) if $K^{(\ell)}$ is invertible for every $\ell=1, \ldots, L$ (notice that we obviously remove $K^{(0)} = 0$ from the list).

\subsubsection{Functions $\cM$ and $\cV$}

For a neater statement of our main result, it is convenient to introduce the following notation. Given a function $h: \R^S \to \R^T$ with polynomial growth of some order $q \ge 0$, we define two induced functions,
\begin{eqnarray}
 \cM_h : \Sym^S_+ \to \Sym_+^T, \quad \cM_h(A) = \EE\sqa{ \bra{ h( \cN_{A} )}^{\otimes 2}},\\
 \cV_h : \Sym^S_+ \to \Sym_+^{T\times T}, \quad \cV_h(A) = \Cov\bra{ \bra{ h(\cN_{A})}^{\otimes 2}}.
\end{eqnarray}
 The letters $\cM$ and $\cV$ remind respectively of ``mean'' and ``variance'' (of $(h(\cN_A))^{\otimes 2}$). Using the polynomial growth assumption and the properties of the norm of tensor products,  is not difficult to argue that, for every $A \in \Sym_+^S$,
 \begin{equation}\label{eq:norm-cm-cv}
  \nor{ \cM_h(A)} \le \cost \nor{h}_{q}^2 (1+\nor{A}^{q}), \quad  \nor{ \cV_h(A)} \le \cost \nor{h}_{q}^4 (1+\nor{A}^{2q}),
 \end{equation}
for some constant $\cost =\cost(q, S)< \infty$.

We anticipate here that in \eqref{eq:derivative-cmh} below we argue that $\cM_h$ is differentiable at any $A \in \Sym^S_+$ that is invertible (i.e., strictly positive definite). We denote with
\begin{equation}
 D_{\cM_h}(A): \Sym^S \to \Sym^T, \quad B \mapsto D_{\cM_h}(A)(B)
\end{equation}
its differential, which is a bounded linear operator such that for every $B \in \Sym^S$,
\begin{equation}\label{eq:norm-d-cm-k}
 \nor{ D_{\cM_h}(A) B} \le \cost \nor{h}_{q}^2 \nor{B},
\end{equation}
with $\cost = \cost(q, S, A)< \infty$ (depending in particular on the smallest eigenvalue of $A$).

Finally, let us notice that we can rewrite the  recursion \eqref{eq:k-ell} defining the infinite width covariance as follows (if $\ell >1$):
\begin{equation}
 K^{(\ell)} = \tr_{a_\ell} \bra{ \cM_{\sigma^{(\ell)}} (K^{(\ell-1)}) }.
\end{equation}

\section{Technical results}\label{sec:technical}

This section is devoted to establish the majority of the technical results utilized in the proofs of our main theorems. For readers who prefer a more concise understanding, we recommend skipping  it entirely, or alternatively starting with  \Cref{sec:proof} and \Cref{sec:bayes}.

\subsection{Quantitative $\delta$-method}

The  $\delta$-method is a technique in probability theory and statistics to describe the law  of a smooth function of approximately Gaussian random variables. In this section, we establish some quantitative versions in terms of the Wasserstein distance for specific functions that we employ in the proof of our main result, \Cref{thm:main-induction}.

\subsubsection{$\delta$-method for the matrix square root.}
 \begin{proposition}\label{prop:delta-method-sqrt}
   Let $A_0 \in \Sym_+^S$ be invertible, let $\cN_{V_0}$ be a Gaussian random variable with values in $\Sym^{S}$ and $A_1$ be a random variable with values in $\Sym^{S}_+$. Then,
   \begin{equation}
    \W_p\bra{ \sqrt{A_1}, \sqrt{ A_0} + D_{\sqrt{\cdot}}(A_0) \cN_{V_0}  } \le \cost \bra{  \nor{ A_1 -A_0}_{L^{2p}}^2 + \W_p\bra{ A_1,  A_0 +\cN_{V_0} } }.
   \end{equation}
   where $\cost = \cost(p, S, A_0)< \infty$.
 \end{proposition}
 \begin{proof}
 On the event
\begin{equation}
 E := \cur{ \nor{A_1 - A_0} \le \lambda(A_0)/2}.
\end{equation}
we apply \eqref{eq:taylor-sqrt} to obtain
\begin{equation}
 \nor{ \sqrt{A_1} - \sqrt{A_0} - D_{\sqrt{\cdot}}(A_0) (A_1 -A_0) } \le  \cost \nor{ A_1 -A_0}^2.
\end{equation}
On the event $E^c$ we apply instead \eqref{eq:ando} and \eqref{eq:norm-dsqrta}, to obtain
\begin{equation}
  \nor{ \sqrt{A_1} - \sqrt{A_0} - D_{\sqrt{\cdot}}(A_0) (A_1 -A_0) } \le \cost \nor{A_1 - A_0}.
\end{equation}
Using Markov's inequality, we also have
\begin{equation}
 \mathbb{P}(E^c) \le \cost \nor{ A_1-A_0}_{L^{2p}}^{2p},
\end{equation}
hence combining these inequalities we find
\begin{equation}\begin{split}
 & \nor{ \sqrt{A_1} - \sqrt{A_0} - D_{\sqrt{\cdot}}(A_0) (A_1 -A_0) }_{L^p} \le \\
 & \le \nor{  I_E \bra{ \sqrt{A_1} - \sqrt{A_0} - D_{\sqrt{\cdot}}(A_0) (A_1 -A_0) }}_{L^p} \\
 & \quad + \nor{  I_{E^c} \bra{ \sqrt{A_1} - \sqrt{A_0} - D_{\sqrt{\cdot}}(A_0) (A_1 -A_0) }}_{L^p}\\
 & \le \cost \nor{A_1-A_0}_{L^{2p}}^2 +  \mathbb{P}(E^c)^{1/(2p)} \nor{A_1 - A_0}_{L^{2p}}\\
 & \le \cost \nor{A_1-A_0}_{L^{2p}}^2,
 \end{split}
\end{equation}
 Consider then any joint realization of $A_1$ and $\cN_{V_0}$. By the triangle inequality,
\begin{equation}
\begin{split}
& \W_p\bra{ \sqrt{A_1},\sqrt{A_0} + D_{\sqrt{\cdot}}(A_0) \cN_{V_0} }  \le  \nor{ \sqrt{A_1} - \sqrt{A_0} - D_{\sqrt{\cdot}}(A_0) \cN_{V_0} }_{L^p}
\\
& \le \nor{  \sqrt{A_1} - \sqrt{A_0} - D_{\sqrt{\cdot}}(A_0) (A_1 -A_0) + D_{\sqrt{\cdot}}(A_0) \bra{ A_1 - A_0 - \cN_{V_0} }}_{L^p}\\
& \le \nor{  \sqrt{A_1} - \sqrt{A_0} - D_{\sqrt{\cdot}}(A_0) (A_1 -A_0)}_{L^p } + \nor{ D_{\sqrt{\cdot}}(A_0) \bra{ A_1 - A_0 - \cN_{V_0} }}_{L^p}\\
 &  \le \cost  \bra{ \nor{ A_1 -A_0}_{L^{2p}}^2 + \nor{A_1 - A_0 - \cN_{V_0} }_{L^p}}
 \end{split}
\end{equation}
having used \eqref{eq:norm-dsqrta}. The thesis follows by minimization upon the joint realizations of $A_1$ and $\cN_{V_0}$.
 \end{proof}

\subsubsection{$\delta$-method for $\cM_h$}

 In this section we first collect some general  facts on the maps $\cM_h$, $\cV_h$ and related quantities and then  we show some results akin to \Cref{prop:delta-method-sqrt}.

 Let us consider a general function $h: \R^S \to \R^T$ (for finite sets $S$, $T$) with polynomial growth of order $q \ge 0$. For $A \in \Sym^S_+$, consider the random variable $h(\cN_A)\stackrel{law}{=} h(\sqrt{A} \cN_S)$. Using \eqref{eq:linear-growth}, we obtain immediately that, for every $p \ge 1$, it holds
\begin{equation}
 \nor{ h(\cN_A)}_{L^p }\le \cost (1 + \nor{\cN_A}_{L^{qp}}^{q} )\le \cost  \bra{1 + \nor{A}^{q/2}},
\end{equation}
with $\cost = \cost(p, q, S)<\infty$. Combining this bound with the identity $\nor{ x^{\otimes k}} = \nor{x}^k$, we obtain that
\begin{equation}\label{eq:nor-lp-hna}
 \nor{ h(\cN_A)^{\otimes k}} _{L^p} \le \cost  \bra{ 1 + \nor{A}^{qk/2}},
\end{equation}
with $\cost = \cost(p,q,k, S)< \infty$. Choosing in particular $p=1$, $k=2$ yield
\begin{equation}\label{eq:nor-cma}
 \nor{ \cM_h(A)} \le \nor{ h(\cN_A)^{\otimes 2}} _{L^1} \le  \cost  \bra{ 1 + \nor{A}^q},
\end{equation}
and choosing $p=1$, $k=4$ and $p=2$, $k=2$ yields
\begin{equation}\label{eq:nor-cva}
 \nor{ \cV_h(A)} \le \nor{ h(\cN_A)^{\otimes 4}} _{L^1} + \nor{ \cM_h(A)}^2  \le  \cost  \bra{ 1 + \nor{A}^{2q}}.
\end{equation}

Next, we investigate the regularity properties of $\cM_h$ and $\cV_h$. We begin with a general differentiability result that we then specialize in our setting.

\begin{lemma}\label{lem:ibp}
 Let $f: \R^S \to \R^T$ have polynomial growth of order $q\ge 0$ and let $A_0 \in \Sym_+^S$ be invertible. Then, for every $B \in \Sym^S$ with $\nor{B} \le \lambda(A_0)/2$, it holds
 \begin{equation}\label{eq:taylor-f-ibp-general}
  \nor{ \EE\sqa{ f( \cN_{A_0+ B})} - \EE\sqa{f(\cN_{A_0})} - D_{\tilde{f}}(A_0)B } \le \cost \nor{f}_q \nor{B}^2,
 \end{equation}
where  $Sym^S \ni B \mapsto D_{\tilde f}(A_0) B \in \Sym^T$ is linear  with
 \begin{equation}\label{eq:derivative-tilde-f}
  \nor{ D_{\tilde f}(A_0) B } \le \cost \nor{f}_q \nor{B},
 \end{equation}
  and $\cost = \cost(q, S, A_0) < \infty$.
\end{lemma}

\begin{remark}
 In particular, we also deduce that
 \begin{equation}\label{eq:f-lipschitz}
  \nor{ \EE\sqa{ f( \cN_{A_0+ B})} - \EE\sqa{f(\cN_{A_0})}} \le \cost \nor{f}_q \nor{B}.
 \end{equation}
\end{remark}

\begin{proof}
The assumption $\nor{B} \le \lambda(A_0)/2 $ ensures that $A_0 + B \in \Sym_+^S$ is invertible, hence the law of $\cN_{A_0 + B}$ is absolutely continuous with respect to that of $\cN_{A_0}$ with density given by
\begin{equation}
\R^S \ni v \mapsto g(v):= \exp\bra{ - \frac 1 2\ang{v, \bra{ (A_0+B)^{-1} - A_0^{-1}} v}} \sqrt{ \frac{\det A_0} {\det (A_0+B)}}.
\end{equation}
Hence, we have the identity
\begin{equation}
\EE\sqa{ f(\cN_{A_0 + B})}  = \EE\sqa{ f( \cN_{A_0}) g(\cN_{A_0})} = \EE\sqa{ f(\sqrt{A_0} \cN ) g(\sqrt{A_0} \cN)},
\end{equation}
where $\cN = \cN_S$ denotes a random variable with standard Gaussian law on $\R^S$. Notice that
\begin{equation}\begin{split}
 g ( \sqrt{A_0} v ) &=  \exp\bra{ - \frac 1 2\ang{ \sqrt{A_0} v, \bra{ (A_0+B)^{-1} - A_0^{-1}} \sqrt{A_0} v}} \sqrt{ \frac{\det A_0} {\det (A_0+B)}}\\
 & = \exp\bra{ - \frac 1 2\ang{ v, \bra{ (\Id_S+B_0)^{-1} - \Id_S}  v}} \sqrt{ \frac{1} {\det (\Id_S +B_0)}},
 \end{split}
\end{equation}
where $B_0 := A_0^{-1/2} B A_0^{-1/2} \in \Sym^S_+$ is such that
\begin{equation}\label{eq:nor-b0}
 \nor{B_0} \le \frac{\nor{B}}{\lambda(A_0)} \le \frac 1 2.
\end{equation}
Thus, if we define for $t \in [0,1]$ and $v \in \R^S$ the function
\begin{equation}
 g(t, v) := \exp\bra{ - \frac 1 2\ang{ v, \bra{ (\Id_S+tB_0)^{-1} - \Id_S}  v}}  \frac{1} {\sqrt{\det (\Id_S +tB_0)}},
\end{equation}
we can write
\begin{equation}
 \EE\sqa{ f(\cN_{A_0 + B})} - \EE\sqa{ f(\cN_{A_0})} = \EE\sqa{ f(\sqrt{A_0} \cN) g(1, \cN)} - \EE\sqa{ f(\sqrt{A_0} \cN) g(0, \cN)}.
\end{equation}
The proof of \eqref{eq:taylor-f-ibp-general} and \eqref{eq:derivative-tilde-f} follows  by differentiating $t \mapsto  \EE\sqa{ f(\sqrt{A_0} \cN) g(t, \cN)}$. Indeed, $t \mapsto g(t, v)$ is smooth, with
\begin{equation}\begin{split}
 \frac{ \partial  g}{\partial t}(t,v)& = -\frac 1 2 \sqa{ \ang{v, (\Id_S + t B_0)^{-2} B_0 v } + \tr\bra{ (\Id_S + t B_0)^{-1} B_0}} g(t,v)\\
 & := m(t,v) g(t,v),
 \end{split}
\end{equation}
which evaluated at $t = 0$ gives
\begin{equation}
  \frac{ \partial  g}{\partial t}(0,v) = m(0,v) =   -\frac 1 2 \sqa{ \ang{v, B_0 v } + \tr(B_0) },
\end{equation}
and is easily bounded from above:
\begin{equation}\label{eq:bound-m0}
 \abs{m(0,v)} \le \cost \nor{v}^2 \nor{B},
\end{equation}
where $\cost  = \cost(S,A_0)<\infty$. We further see that
\begin{equation}
  \frac{ \partial^2  g}{\partial t^2}(t,v) =  \sqa{ m(t,v)^2 + \frac{ \partial  m}{\partial t}(t,v) } g(t,v),
\end{equation}
and one easily realizes (using also the condition \eqref{eq:nor-b0}) that, for some constant $\cost = \cost(S, A_0)< \infty$ it holds
\begin{equation}
 \sup_{t \in [0,1]}  \abs{ m(t,v)^2 + \frac{ \partial  m}{\partial t}(t,v) } \le \cost  \nor{v}^4 \nor{B_0}^2 \le \cost \nor{v}^4 \nor{B}^2.
\end{equation}
We obtain therefore that
\begin{equation}
\begin{split}
 \nor{ \EE\sqa{ f(\sqrt{A_0} \cN) \int_0^1  \frac{ \partial^2  g}{\partial t^2}(t, \cN) (1-t) dt } } & \le \EE\sqa{ \nor{f}_{q} \bra{ 1+ \nor{\cN_{A_0 + tB}}^q}   \nor{\cN_{A_0 + tB}}^4 \nor{B}^2} \\
 &  \le \cost \nor{f}_q \nor{B}^2,
 \end{split}
\end{equation}
where $\cost = \cost(q,S, A_0)< \infty$, hence by Taylor's formula, we get the bound
\begin{equation}
 \begin{split}
 &  \nor{ \EE\sqa{ f(\sqrt{A_0} \cN) g(1, \cN)} - \EE\sqa{ f(\sqrt{A_0} \cN) g(0, \cN)} - \EE\sqa{ f(\sqrt{A_0}) \frac{ \partial  g}{\partial t}(0,v) } } \le \\
  & \le  \nor{ \EE\sqa{ f(\sqrt{A_0} \cN) \int_0^1  \frac{ \partial^2  g}{\partial t^2}(t,\cN) (1-t) dt } } \\
  & \le \cost \nor{f}_q \nor{B}^2.
 \end{split}
\end{equation}
Thus, to conclude it is sufficient to define
\begin{equation}
 D_{\tilde f }(A_0)(B) :=- \frac 1 2  \EE\sqa{ f(\sqrt{A_0} \cN) \bra{ \ang{\cN, A_0^{-1/2} B A_0^{-1/2} \cN } + \tr(A_0^{-1/2}B A_0^{-1/2}) }}
\end{equation}
and finally check that inequality \eqref{eq:derivative-tilde-f} holds, using \eqref{eq:bound-m0}.
\end{proof}

By specializing the above result to the cases $f = h^{\otimes 2}$ with $h: \R^S \to \R^T$ with polynomial growth of order $q\ge 0$, we obtain that, for every invertible $A_0 \in \Sym_+^S$ and $B \in \Sym^S$ with $\nor{B}_{op}\le \lambda(A_0)/2$, it holds
\begin{equation}\label{eq:derivative-cmh}
 \nor{  \cM_h(A_0 +B) - \cM_h(A_0) - D_{\cM_h}(A_0) B } \le \cost \nor{h}_q^{2} \nor{B}^2,
\end{equation}
where for every $B \in \Sym^S$, $B \mapsto D_{\cM_h}(A_0) = D_{\widetilde{h^{\otimes 2}} }(A_0)B$ is linear with
\begin{equation}\label{eq:nor-derivative-cmh}
  \nor{ D_{\cM_h}(A_0) B} \le \cost \nor{h}_q^2\nor{B},
\end{equation}
and $\cost = \cost(q, S, A_0)<\infty$.

We combine the above differentiability result with the simpler bound \eqref{eq:nor-cma} to establish the following quantitative $\delta$-method for the map  $\cM_h$.

\begin{proposition}\label{prop:delta-method-mh}
   Let $A_0 \in \Sym_+^S$ be invertible, let $\cN_{V_0}$ be a Gaussian random variable with values in $\Sym^{S}$ and $A_1$ be a random variable with values in $\Sym^{S}_+$. Then, for every function $h: \R^S \to \R^T$ with polynomial growth of order $q\ge 2$, it holds
   \begin{equation}\begin{split}
    & \W_p\bra{ \cM_h(A_1), \cM_h(A_0) + D_{\cM_h}(A_0) \cN_{V_0}  } \le\\
    & \quad  \le \cost \nor{h}_{q}^2  \bra{  \nor{ A_1 -A_0}_{L^{2p}}^2 + \nor{ A_1 -A_0}_{L^{qp}}^q + \W_p\bra{ A_1,  A_0 +\cN_{V_0} } }.
    \end{split}
   \end{equation}
   where $\cost = \cost(p,q,A_0)< \infty$.
 \end{proposition}

 The condition $q\ge 2$ can be obviously relaxed, since any function with polynomial growth of order $q\le 2$ has also polynomial growth of order $2$.

\begin{proof}
 We assume for simplicity that $\nor{h}_q \le 1$. We introduce the event
 \begin{equation}
  E = \cur{ \nor{A_1-A_0} \le \lambda(A_0)/2}.
\end{equation}
If $E$ holds, we can apply \eqref{eq:derivative-cmh} with $B:= A_1-A_0$. Otherwise, we apply \eqref{eq:nor-cma} with $A:=A_1$ and $A:=A_0$ and also \eqref{eq:nor-derivative-cmh} with $B:=A_1-A_0$, yielding
\begin{equation}
\nor{ \cM_h(A_1) - \cM_h(A_0) -D_{\cM_h}(A_0) (A_1-A_0)} \le \cost (1+ \nor{A_1-A_0}^q).
\end{equation}
By Markov's inequality, we have $\mathbb{P}(E^c) \le \cost \nor{A_1-A_0}_{L^{{2p}}}^{{2p}}$. Thus,
\begin{equation}\label{eq:split-markov}\begin{split}
 & \nor{ \cM_h(A_1) - \cM_h(A_0) -D_{\cM_h}(A_0) (A_1-A_0)}_{L^p} \le \\
 & \le \nor{ I_E\bra{ \cM_h(A_1) - \cM_h(A_0) -D_{\cM_h}(A_0) (A_1-A_0)} }_{L^p} + \\
  & \quad + \nor{ I_{E^c}\bra{ \cM_h(A_1) - \cM_h(A_0) -D_{\cM_h}(A_0) (A_1-A_0)} }_{L^p}\\
  & \le \cost \bra{\nor{A_1-A_0}_{L^{2p}}^2 + \nor{A_1-A_0}_{2p}^2 + \nor{A_1-A_0}_{qp}^q}\\
  & \le \cost \bra{ \nor{A_1-A_0}_{L^{2p}}^2 +\nor{A_1-A_0}_{qp}^q }
 \end{split}
\end{equation}
Next, we consider any joint realization of $A_1$ and $\cN_{V_0}$ and bound from above
\begin{equation}
 \begin{split}
  & \W_p\bra{ \cM_h(A_1), \cM_h(A_0) + D_{\cM_h}(A_0) \cN_{V_0} }  \le\\
  & \le \nor{ \cM_h(A_1)- \cM_h(A_0) - D_{\cM_h}(A_0) \cN_{V_0}}_{L^p} \\
  & \le \nor{\cM_h(A_1)- \cM_h(A_0) - D_{\cM_h}(A_0) (A_1-A_0) }_{L^p} + \nor{D_{\cM_h}(A_0) \bra{ A_1 -A_0 -\cN_{V_0}} }_{L^p}\\
  & \le \cost \bra{ \nor{A_1-A_0}_{L^{2p}}^2 +\nor{A_1-A_0}_{qp}^q  + \nor{ A_1 -A_0 -\cN_{V_0}}_{L^p}},
 \end{split}
\end{equation}
having also used \eqref{eq:nor-derivative-cmh} with $B:=A_1 - A_0 - \cN_{V_0}$. Minimization upon the realizations of $A_1$ and $\cN_{V_0}$ yields the thesis.
\end{proof}

We next focus on the function $\cV_h$. Indeed, it is straightforward to prove
\begin{equation}\label{eq:cv-lipschitz}
 \nor{ \cV_h(A_0 + B ) - \cV_h(A_0) } \le \cost \nor{h}^4_q \nor{B},
\end{equation}
where $\cost = \cost(q,S,A_0)< \infty$, assuming that $A_0$ is invertible and $\nor{B} \le \lambda(A)/2$. Indeed, by \eqref{eq:f-lipschitz} with $f = h^{\otimes 4}$, we obtain
\begin{equation}
\nor{ \EE\sqa{ ( h(A_0 + B ) )^{\otimes 4} } -  \EE\sqa{ ( h(A_0) )^{\otimes 4} } } \le \cost \nor{h}_{q}^4 \nor{B},
\end{equation}
using the fact that $h^{\otimes 4}$ has polynomial growth of order $4q$. Then, combining \eqref{eq:nor-cma} with \eqref{eq:f-lipschitz} applied with $f = h^{\otimes 2}$ we obtain similarly that
\begin{equation}
\nor{ \cM_h(A_0 +B)^{\otimes 2} - \cM_h(A_0)^{\otimes 2}} \le \cost \nor{h}_q^4 \nor{B},
\end{equation}
hence using these two bounds we get  \eqref{eq:cv-lipschitz}. In the following lemma we combine such inequality with \eqref{eq:ando} to obtain a similar bound for $\sqrt{\cV_h}$.

\begin{lemma}
 Let $A_0 \in \Sym_+^S$ be invertible and let $A_1$ be a random variable with values in $\Sym_+^S$. Then, for every $h: \R^S \to \R^T$ with polynomial growth of order $q\ge 1$, it holds
 \begin{equation}\label{eq:lipschitz-sqrt-cvh}
  \nor{ \sqrt{\cV_h(A_1)} - \sqrt{\cV_h(A_0)}}_{L^p} \le \cost \bra{ \nor{ A_1 -A_0}_{L^p} + \nor{A_1 - A_0}_{L^{2qp}}^{2q}}
\end{equation}
where $\cost = \cost(p,q, S, T, A_0, h)< \infty$.
\end{lemma}

\begin{proof}
Let us notice first the following fact. If $f: \R^S \to \R^T$ is Borel (with polynomial growth) and $A_1$, $A_0$ are both invertible, then $\Cov ( f(\cN_{A_1})) \ll \Cov(f(\cN_{A_0}))$. Indeed, if $v \in \R^{T}$ is such that  $\Cov( f(\cN_{A_0})) v = 0$, then also  the variance $\operatorname{Var}( \ang{ f(\cN_{A_0}), v}) = \ang{ v, \Cov( f(\cN_{A_0})) v }  = 0$, i.e., the random variable $\ang{f(\cN_{A_0}), v}$ is a.s.\ constant, or equivalently the function $\R^S \ni w \mapsto \ang{f(w) v}$ is almost everywhere constant with respect to the (centered) Gaussian measure with covariance $A_0$. But such measure is equivalent to the Lebesgue measure, hence also equivalent to the centered Gaussian measure with covariance $A_1$. Thus, we have  that  also $\ang{f(\cN_{A_1}), v}$ is a.s.\ constant, hence  $\operatorname{Var}( \ang{f(\cN_{A_1}), v}) = \ang{v, \Cov( f(\cN_{A_1})) v }  = 0$ and therefore (since the covariance is positive semidefinite) $\Cov( f(\cN_{A_1})) v = 0$.

Back to the proof of \eqref{eq:lipschitz-sqrt-cvh}, we assume without loss of generality that $\nor{h}_q \le 1$ and introduce the event
\begin{equation}
  E = \cur{ \nor{A_1 - A_0} \le \lambda(A_0)/2}
\end{equation}
so that, if $E$ holds then $A_1$ is also invertible and therefore, by the above considerations applied to $f = h^{\otimes 2}$, it holds $\cV_h(A_1) \ll \cV_h(A_0)$. If $E$ holds, we can thus apply \eqref{eq:ando} with $B:=\cV_h(A_1)$ and $A := \cV_h(A_0)$ to obtain
\begin{equation}
 \nor{ \sqrt{\cV_h(A_1)} - \sqrt{\cV_h(A_0)}} \le \cost \nor{ \cV_h(A_1) - \cV_h(A_0)}.
\end{equation}
where $\cost = \cost \bra{\lambda_+\bra{\cV_h(A_0)}}< \infty$ (the case $\cV_h(A_0) = 0$ is trivial, for then also $\cV_h(A_1) = 0$). Using \eqref{eq:cv-lipschitz} with $A:=A_1$ and the fact that $\nor{A_1-A_0} \le \lambda(A_0)/2 \le \cost$ on $E$, we next find
\begin{equation}
 \nor{ \sqrt{\cV_h(A_1)} - \sqrt{\cV_h(A_0)}} \le \cost \nor{ A_1 -A_0}.
\end{equation}
On $E^c$, we simply use \eqref{eq:holder} and then \eqref{eq:nor-cva}, so that
\begin{equation}\begin{split}
 \nor{ \sqrt{\cV_h(A_1)} - \sqrt{\cV_h(A_0)}} & \le \cost \sqrt{ \nor{\cV_h(A_1)} +\nor{\cV_h(A_0) } } \\
 & \le \cost  (1 + \nor{A_1 - A_0}^{2q} ).
 \end{split}
\end{equation}
Therefore, splitting the integration over $E$ and $E^c$ and using Markov's inequality as in \eqref{eq:split-markov} (with exponent $1$ instead of $2$ when integrating over $E$ and exponent $2q$ instead of $q$ on $E^c$) we find the thesis.
\end{proof}

 Ultimately, we are going to impose that $h$ is Lipschitz, however it seemed more natural so far to establish the bounds in this section using only the growth of $h$, and not its regularity. Of course, if $h$ is Lipschitz, then one can straightforwardly obtain stronger bounds. For example,
 we have the inequality
 \begin{equation}
  \nor{ h^{\otimes 2} (\sqrt{A_1} \cN_S) -  h^{\otimes 2} (\sqrt{A_0} \cN_S)}_{L^p} \le \cost \nor{h}_{\Lip}^2 \bra{ \nor{\sqrt{A_1}} +\nor{\sqrt{A_0}} } \nor{ \sqrt{A_1} - \sqrt{A_0}}
 \end{equation}
 where $\cost = \cost(p,S)< \infty$, and combining this inequality with \eqref{eq:ando} (e.g.\ if $A_0$ is invertible) we find
 \begin{equation}\label{eq:lipschitz-h-lp}
  \nor{ h^{\otimes 2} (\sqrt{A_1} \cN_S) -  h^{\otimes 2} (\sqrt{A_0} \cN_S)}_{L^p} \le \cost (1 + \sqrt{ \nor{A_1 - A_0} }) \nor{A_1 - A_0},
 \end{equation}
 where $\cost = \cost(p,S, A_0, h) < \infty$.

 We end this section by giving the simple proof of the following bound, used in the proof of \eqref{eq:induction-kernel-basteri}: given random variables $X$, $Y$ taking values in $\R^{n \times S}$, for every Lipschitz continuous $h: \R^S \to \R^T$ it holds
 \begin{equation}\label{eq:kernel-basteri-tool}
  \W_p\bra{ \sum_{i=1}^n h(X[i,:])^{\otimes 2}, \sum_{i=1}^n h(Y[i,:])^{\otimes 2}} \le \cost \nor{h}_{\Lip}^2 \bra{ 1 + \nor{Y}_{L^{2p}} + \W_{2p}\bra{X,Y}} \W_{2p}\bra{X,Y},
 \end{equation}
 where $\cost = \cost(p) < \infty$.  Indeed, given any joint realization of $X$ and $Y$, by adding and subtracting the quantity $\sum_{i=1}^n  h(X[i,:]) \otimes h(Y[i,:])$, we find
 \begin{equation}\begin{split}
                  &  \nor{ \sum_{i=1}^n h(X[i,:])^{\otimes 2}-  \sum_{i=1}^n h(Y[i,:])^{\otimes 2} }_{L^p } \le\\
                  & \le \nor{ \sum_{i=1}^n \nor{h(X[i,:])}\nor{ h(X[i,:])-h(Y[i,:])}} _{L^p} \\
                  & \quad + \nor{ \sum_{i=1}^n \nor{h(Y[i,:])}\nor{ h(X[i,:])-h(Y[i,:])}} _{L^p}   \\
                  & \le  \nor{  \nor{h(X)} \nor{ h(X)-h(Y)}} _{L^p} + \nor{  \nor{h(Y)} \nor{ h(X)-h(Y)}} _{L^p},
                 \end{split}
 \end{equation}
having applied Cauchy-Schwarz inequality and writing $h(X) := (h(X[i,:]))_{i=1, \ldots, n}$ and similarly for $h(Y)$. Focusing on the first term (the second is similar), we find
\begin{equation}\begin{split}
  \nor{  \nor{h(X)} \nor{ h(X)-h(Y)}} _{L^p} & \le \nor{h(X)}_{L^{2p}}  \nor{ h(X)-h(Y)}_{L^{2p}} \\
  & \le \nor{h}_{\Lip}^2 (1 + \nor{X}_{L^{2p}}) \nor{X-Y}_{L^p},
  \end{split}
\end{equation}
having used the fact that $h$ is Lipschitz continuous and acts componentwise in the definition of $h(X)$. Finally, we bound
\begin{equation}
 \nor{X}_{L^{2p}} \le \nor{Y}_{L^{2p}} + \nor{X-Y}_{L^{2p}}
\end{equation}
and \eqref{eq:kernel-basteri-tool} follows by minimization over the joint definitions of $X$ and $Y$.

\subsection{Quantitative central limit theorem}

In this section we establish the following quantitative version of the central limit theorem for a ``sum of squares'' of Lipschitz functions of conditionally independent Gaussian variables.

\begin{proposition}\label{prop:clt-new}
Let $n \ge 1$, $S, T$ be  finite sets, let $A_0 \in \Sym_+^S$ be invertible and $A_1$  be a random variable with values in $\Sym^S_+$. Let $\cN_{n\times S}$ denote a standard Gaussian random variable with values in $\R^{n\times S}$ and independent of $A_1$ and let $h: \R^S \to \R^T$  be Lipschitz continuous. Then,
\begin{equation}\label{eq:clt-quant-square}
  \W_p\bra{ \frac 1 n \sum_{i=1}^n h( \sqrt{A_1} \cN_{n\times S}[i,:])^{\otimes 2 }, \cM_h(A_1) + \frac{\cN_{\cV_h(A_1)}}{\sqrt{n}} } \le \cost \bra{ \frac{1}{n} +  \nor{A_1 -A_0}_{L^{2p}}^2 },
\end{equation}
where $\cost = \cost(p, S, T, A_0, h) < \infty$.
\end{proposition}

Before we address the proof, let us first recall and combine together two general results that we employ in our argument. Consider a random variable $Z$ with values in $\R^d$,  let $(Z_i)_{i=1}^n$ be independent random variables, each with the same law as $Z$, and write
\begin{equation}
 \bar{Z}_n :=\frac 1 n \sum_{i=1}^n Z_i
\end{equation}
for their mean. The first result is the classical Rosenthal's inequality (see e.g.\ \cite{osekowski2012note}) which yields for  any $p\ge 2$ the following bound:
\begin{equation}\label{eq:rosenthal}
 \nor{ \bar{Z}_n - \EE\sqa{Z}}_{L^p} \le \cost \bra{  \frac{ \nor{Z-  \EE\sqa{Z}}_{L^2} }{\sqrt{n}} + \frac{ \nor{Z- \EE\sqa{Z}}_{L^p} }{n^{1-1/p}} }
 \le \cost \frac{ \nor{Z- \EE\sqa{Z}}_{L^p} }{\sqrt{n}},
\end{equation}
where $\cost = \cost(p)<\infty$. The second result that we recall is a the quantitative central limit theorem with respect to the Wasserstein distance, precisely \cite[Theorem 1]{bonis_steins_2020}, which in our notation  and after some straightforward reductions reads as follows, again for any $p \ge 2$:
\begin{equation}\label{eq:bonis}
 \W_p \bra{ \bar{Z}_n,  \EE\sqa{Z} + \frac{\cN_d}{\sqrt{n} }} \le \cost \frac{ \nor{Z}_{L^{p+2}}^{1+2/p}  }{n},
\end{equation}
where $\cost =\cost(p,d)< \infty$, and provided that $\Cov(Z) = \Id_d$. Notice for example that for $p=2$ it requires a finite absolute moment of order $4$ to be effective.

In our application to the proof of \Cref{prop:clt-new}, we need however to slightly extend  \eqref{eq:bonis} allowing for general covariances (possibly not invertible): this is achieved by decomposing the space $\R^d$ into the orthogonal sum of two subspaces and respectively apply \eqref{eq:rosenthal} and \eqref{eq:bonis} on each subspace. For simplicity, and without loss of generality, we assume that $Z$ (hence each random variable $Z_i$) is centered.  Given $d_0 \in \cur{0,1,\ldots, d}$, we consider the decomposition $\R^d = \R^{d_0} \times \R^{d - d_0}$ and write  accordingly
\begin{equation}
 Z = (Z^{0}, Z^{\perp}),  \quad Z_i = (Z^0_i, Z_i^{\perp}) \quad \bar{Z}_n = ( \bar{Z}^0_n, \bar{Z}^{\perp}_n)
\end{equation}
 and similarly $\cN_{\Cov(Z)} = (\cN_{\Cov(Z)}^0, \cN_{\Cov(Z)}^{\perp})$. Notice that both $\cN_{\Cov(Z)}^0$ and $\cN_{\Cov(Z)}^{\perp}$ are Gaussian variables (in general they are not independent) with covariances given by the respective top $d_0 \times d_0$ and bottom $(d-d_0) \times (d-d_0)$ diagonal blocks of $\Cov(Z)$, which coincide with $\Cov(Z^0)$ and $\Cov(Z^\perp)$ respectively.

Using the triangle inequality, we then bound from above
\begin{equation}\begin{split}
  \W_p\bra{ \bar{Z}_n,   \frac{\cN_{\Cov(Z)}}{\sqrt{n} } } &
 \le \W_p\bra{ \bar Z _n,   \bra{ \bar{Z}^0_n, 0} }  \\
& \quad + \W_p\bra{  \bra{ \bar{Z}_n^0,  0}, \bra{ \frac{\cN_{\Cov(Z)}^0} {\sqrt{n} }, 0 }} \\
 & \quad + \W_p\bra{   \bra{ \frac{ \cN_{\Cov(Z)}^0}{\sqrt{n} }, 0}, \frac{\cN_{\Cov(Z)}}{\sqrt{n} } },
\end{split}
\end{equation}
and we estimate separately the three terms. For the last one, we simply use \eqref{eq:w-2-trivial-bound} so that
\begin{equation}
\begin{split}
 \W_p\bra{  \bra{ \frac{\cN_{\Cov(Z)}^0}{\sqrt{n} }, 0}, \frac{\cN_{\Cov(Z)}}{\sqrt{n} }  }  &\le \frac 1 {\sqrt{n}} \W_p\bra{  \bra{ \cN_{\Cov(Z)}^0, 0},\bra{ \cN_{\Cov(Z)}^0, \cN_{\Cov(Z)}^\perp}  } \\
 & \le \frac{1}{\sqrt{n}} \nor{ \cN_{\Cov(Z)}^{\perp} }_{L^p} \le \cost  \frac{ \nor{ \sqrt{ \Cov(Z^\perp) }}}{\sqrt{n}}\\
 & \le \cost \frac{\nor{Z^\perp}_{L^2 }}{\sqrt{n}}  \le  \cost  \frac{\nor{Z^\perp}_{L^p }}{\sqrt{n}},
 \end{split}
\end{equation}
where $\cost = \cost(p)< \infty$. We estimate similarly the first term, using again \eqref{eq:w-2-trivial-bound} but then Rosenthal's inequality \eqref{eq:rosenthal}, obtaining
\begin{equation}\begin{split}
 \W_p\bra{ \bar Z_n,   \bra{ \bar{Z}_n^0, 0} }  & = \W_p\bra{  \bra{\bar{Z}_n^0, \bar{Z}_n^\perp } ,  \bra{ \bar{Z}_n^0, 0} }\\
 & \le \nor{\bar{Z}_n^{\perp} }_{L^p} \le \cost  \frac{\nor{Z^\perp}_{L^p }}{\sqrt{n}},
 \end{split}
\end{equation}
with $\cost = \cost(p)< \infty$. For the second term, i.e.,
\begin{equation}
 \W_p\bra{  \bra{ \bar{Z}_n^0,  0}, \bra{\frac{\cN_{\Cov(Z)}^0} {\sqrt{n} }, 0 }} = \W_p\bra{  \bar{Z}_n^0, \frac{\cN_{\Cov(Z)}^0} {\sqrt{n} } }
\end{equation}
we use instead \eqref{eq:bonis}. Precisely, we assume that $\Cov(Z^0)$ is invertible and apply \eqref{eq:bonis} to the standardized variables $\Cov(Z^0)^{-1/2} Z^0_i$ (whose covariance is $\Id_{d_0}$). We obtain
\begin{equation}
\begin{split}
 \W_p\bra{  \bar Z_n^0, \frac{\cN_{\Cov(Z)}^0} {\sqrt{n} } } & \le \sqrt{ \nor{\Cov(Z^0) }_{op} } \W_p\bra{  \frac 1 n \sum_{i=1}^n \Cov(Z^0)^{-1/2} Z_i^0, \frac{\cN_{d^0}} {\sqrt{n} } }\\
& \le \cost \sqrt{ \nor{\Cov(Z^0) }_{op} }  \frac{ \nor{ \Cov(Z^0)^{-1/2} Z^0}_{L^{p+2}}^{1+2/p}  }{n}\\
& \le \cost \sqrt{  \frac{  \nor{\Cov(Z^0) }_{op} } {\lambda(\Cov(Z^0))^{1+2/p}} } \frac{ \nor{Z}_{L^{p+2}}^{1+2/p}}{n},
\end{split}
\end{equation}
where $\cost = \cost(p,d)<\infty$, having used also that $\nor{Z^0} \le \nor{Z}$. We obtain therefore the following general inequality
\begin{equation}
 \label{eq:clt-abstract-split}
 \W_p\bra{  \frac 1 n \sum_{i=1}^n Z_i,   \frac{\cN_{\Cov(Z)}}{\sqrt{n} } } \le \cost \bra{ \sqrt{  \frac{  \nor{\Cov(Z^0) }_{op} } {\lambda(\Cov(Z^0))^{1+2/p}} } \frac{ \nor{Z}_{L^{p+2}}^{1+2/p}}{n} +  \frac{\nor{Z^\perp}_{L^p }}{\sqrt{n}}  },
\end{equation}
with $\cost = \cost(p,d)<\infty$, under the assumption $\Cov(Z^0)$ is invertible (and only for simplicity that $Z$ is centered). Notice that,  using the trivial decomposition corresponding to $d_0 = 0$, we obtain
\begin{equation}\label{eq:rosenthal-wp}
 \W_p\bra{  \bar{Z}_n,   \frac{\cN_{\Cov(Z)}}{\sqrt{n} } } \le \cost  \frac{\nor{Z}_{L^p }}{\sqrt{n}}.
\end{equation}
We finally notice that \eqref{eq:clt-abstract-split} holds true for more general decompositions of $\R^d$ into  the orthogonal sum of two subspaces, because the Wasserstein distance and the inequalities  \eqref{eq:rosenthal} and \eqref{eq:bonis} are invariant with respect to orthogonal transformations of $\R^d$.

We can now address the proof of \Cref{prop:clt-new}.

\begin{proof}[Proof of \Cref{prop:clt-new}]
Up to a rescaling of $h$ we can assume that $\nor{h}_{\Lip} \le 1$ and also that $\lambda(A_0) =1$.  Moreover, since we argue for a fixed $h$, we simplify the notation and write, for $A \in \Sym_+^S$, $\cM(A):= \cM_h(A)$ and $\cV(A) := \cV_h(A)$. We also introduce the random variables (taking values in $\Sym_+^T \subseteq \R^{T\times T}$, that we naturally identify with $\R^{d}$ with $d = |T|^2$),
\begin{equation}
 Z(A) := h(\sqrt{A}\cN_S)^{\otimes 2}, \quad Z_i(A) := h(\sqrt{A} \cN_{n\times S}[i,:])^{\otimes 2} \quad \bar{Z}_n(A) := \frac 1 n \sum_{i=1}^n Z_i(A),
\end{equation}
so that $(Z_i(A))_{i=1,\ldots, n}$ are independent copies of $Z(A)$ and we have the identities
\begin{equation}
  \EE\sqa{Z(A)} = \cM(A), \quad \Cov(Z(A)) = \cV(A) \quad \text{and}\quad \bar{Z}_n(A) = \frac 1 n \sum_{i=1}^n h( \sqrt{A_1} \cN_{n\times S}[i,:])^{\otimes 2 }.
\end{equation}
Next, for  $\eps\in (0,1]$ we define the set
\begin{equation}
 E = E(\eps) := \cur{A \in \Sym_+^S \, : \,  \nor{ A - A_0} < \eps}.
\end{equation}
and we claim that, for some suitable choice of $\eps = \eps(p,S, T,  A_0, h)\in (0,1]$, it holds
\begin{equation}
\begin{split}\label{eq:claim-prop-clt-square-E}
& \W_p\bra{ \bar{Z}_n(A),  \cM( A) +  \frac{\cN_{\cV(A)}}{\sqrt{n}} }  \le  \cost \bra{\frac{1}{n} + \nor{A-A_0}^2 }, \quad \text{for every $A \in E$,}
\end{split}
\end{equation}
while \begin{equation}\label{eq:claim-prop-clt-square-E-not}
 \W_p\bra{ \bar{Z}_n(A),  \cM( A) + \frac{\cN_{\cV(A)}}{\sqrt{n}} } \le \cost \frac{1 + \nor{A-A_0}}{\sqrt{n}}, \quad \text{for every $A \notin E$.}
\end{equation}
Let us argue immediately that these two inequalities are sufficient to obtain the thesis. Indeed, by  Markov's inequality, we bound from above
\begin{equation}\label{eq:markov-ec-2}
 \mathbb{P}(E^c )\le \frac{ \nor{ A_1 - A_0}_{L^{2p}}^{2p}}{\eps^{2p}} = \cost \nor{ A_1 - A_0}_{L^{2p}}^{2p},
\end{equation}
with $\cost := \eps^{-2p} < \infty$. Then, by \eqref{eq:convexity} and splitting the integration on $E$ and $E^c$, using \eqref{eq:claim-prop-clt-square-E}, \eqref{eq:claim-prop-clt-square-E-not} and  \eqref{eq:markov-ec-2}, we have
\begin{equation}
\begin{split}
 & \W_p\bra{ \bar{Z}(A_1),  \cM( A_1) +  \frac{\cN_{\cV(A_1)}}{\sqrt{n}} } \\
 & \le \bra{ \int_E  \W_p\bra{ \bar{Z}_n(A),  \cM( A) +  \frac{\cN_{\cV(A)}}{\sqrt{n}} }^p d \mathbb{P}_{A_1}(A)  }^{1/p} \\
  & \quad \qquad + \bra{ \int_{E^c} \W_p\bra{\bar{Z}_n(A), \cM( A) +  \frac{\cN_{\cV(A)}}{\sqrt{n}}}^p  d \mathbb{P}_{A_1}(A)  }^{1/p} \\
   & \le  \cost \bra{ \frac{1}{n}  + \nor{A_1 -A_0}_{L^{2p}}^2 + \frac 1 {\sqrt{n}} \mathbb{P}(E^c)^{1/(2p)} \nor{A-A_0}_{L^{2p}}  } \\
  & \le \cost \bra{ \frac 1 n + \nor{A_1-A_0}_{L^{2p}}^2}.
\end{split}
\end{equation}

First, we show \eqref{eq:claim-prop-clt-square-E}: given $A \in E$, we use \eqref{eq:clt-abstract-split} with the orthogonal decomposition induced by  $\R^d = \operatorname{Im} \Cov(Z(A_0)) \oplus \Ker \Cov(Z(A_0))$ (notice that we employ $A_0$, not $A$). Since $\Cov(Z(A_0)) = \cV(A_0)$, by \eqref{eq:cv-lipschitz}, we deduce that
\begin{equation}\label{eq:cov-za-continuous}
 \nor{ \Cov(Z(A))  - \Cov(Z(A_0))} \le \cost  \nor{A-A_0}
\end{equation}
with $\cost = \cost(S, A_0, h)< \infty$. Moreover, recalling the notation from \eqref{eq:clt-abstract-split},
\begin{equation}
\begin{split}
 \nor{ Z(A)^{\perp} - \EE\sqa{ Z(A)}^{\perp} }_{L^p} &\le \nor{ Z(A)^{\perp} -Z(A_0)^{\perp}}_{L^p} + \nor{ \EE\sqa{Z(A)} - \EE\sqa{Z(A_0)}}\\
 & \le \cost  \bra{ \nor{ Z(A) -Z(A_0)}_{L^p} +  \nor{ Z(A) -Z(A_0)}_{L^1}}\\
 & \le \cost \nor{ Z(A) -Z(A_0)}_{L^p},
 \end{split}
\end{equation}
where $\cost = \cost(p)< \infty$, having used in the first inequality the fact that the random variable  $Z(A_0)^{\perp}$ is constant and equal to $\EE\sqa{ Z(A_0)^{\perp}}$, precisely because $\Cov(Z(A_0)^{\perp}) = 0$. Using \eqref{eq:lipschitz-h-lp} (in particular we rely here on the fact that $A_0$ is invertible the assumption that $h$ is Lipschitz continuous) we deduce that
\begin{equation}
  \nor{ Z(A) -Z(A_0)}_{L^p} \le \cost  \bra{1 + \sqrt{\nor{A-A_0}}}\nor{A-A_0} \le \cost \nor{A-A_0},
\end{equation}
where $\cost = \cost(p, S, A_0, h)<\infty$ and we used that $\nor{A-A_0} \le \eps \le 1$. Therefore, we also find
\begin{equation}
 \nor{ Z(A)^{\perp} - \EE\sqa{ Z(A)}^{\perp} }_{L^p} \le \nor{ Z(A) - \EE\sqa{ Z(A)} }_{L^p}  \le \cost \nor{A-A_0},
\end{equation}
where $\cost = \cost(p, S, A_0)< \infty$.

Now if $\Cov(Z(A_0)) = 0$, then the right hand side \eqref{eq:clt-abstract-split} would reduce to the second term only, for which we already have established a proper upper bound:
\begin{equation}
 \W_p\bra{ \bar{Z}_n(A),  \cM(A) + \cN_{\cV(A)} } \le \cost  \frac{\nor{A-A_0}}{\sqrt{n}} \le \cost \bra{ \frac 1 n + \nor{A-A_0}^2},
\end{equation}
and \eqref{eq:claim-prop-clt-square-E} already follows (we can safely choose $\eps = 1$). Otherwise, in the right hand side of \eqref{eq:clt-abstract-split} there appears also the term
\begin{equation}
  \sqrt{  \frac{  \nor{\Cov(Z(A)^0) }_{op} } {\lambda(\Cov(Z(A)^0))^{1+2/p}} } \frac{ \nor{Z(A)-\EE\sqa{Z(A)} }_{L^{p+2}}^{1+2/p}}{n},
\end{equation}
where $\lambda(\Cov(Z(A)^0)>0$. Using \eqref{eq:nor-lp-hna} (with $q=1$ and $k=2$) we deduce that
\begin{equation}
 \nor{Z(A)-\EE\sqa{Z(A)} }_{L^{p+2}} \le \cost (1+\nor{A}^{p+2} ) \le \cost
\end{equation}
having used that $\nor{A} \le \nor{A_0} + \nor{A-A_0} \le \cost$ on $E$ (recall that $\eps \le 1$). Arguing similarly, and using \eqref{eq:cov-za-continuous}, we find
\begin{equation}\begin{split}
 \nor{\Cov(Z(A)^0) }_{op} & \le  \nor{\Cov(Z(A)) } \\
 & \le \cost \bra{ \nor{\Cov(Z(A_0) } + \nor{\Cov(Z(A)) - \Cov(Z(A_0)}} \le \cost.
 \end{split}
\end{equation}
Thus, to conclude with \eqref{eq:claim-prop-clt-square-E}, is it sufficient to choose $\eps$ sufficiently small so that
\begin{equation}
 \lambda(\Cov(Z(A)^0)) \ge \lambda(\Cov(Z(A_0)^0))/2
\end{equation}
and obtain the upper bound
\begin{equation}
 \sqrt{  \frac{  \nor{\Cov(Z(A)^0) }_{op} } {\lambda(\Cov(Z(A)^0))^{1+2/p}} } \frac{ \nor{Z(A)-\EE\sqa{Z(A)} }_{L^{p+2}}^{1+2/p}}{n} \le \cost \cdot \frac{1}{n}.
\end{equation}
Such a choice for $\eps$ is possible, since again by \eqref{eq:cov-za-continuous}
\begin{equation}
 \nor{ \Cov(Z(A)^0) - \Cov(Z(A_0)^0)}\le \nor{ \Cov(Z(A)) - \Cov(Z(A_0))}  \le \bar \cost \nor{A - A_0},
\end{equation}
where $\bar \cost = \bar{\cost}(S,T, A_0, h)$, hence it is sufficient to choose
\begin{equation}
 \eps =  \frac {\lambda(\Cov(Z(A)^0))}{\bar \cost}.
\end{equation}

Finally, we prove \eqref{eq:claim-prop-clt-square-E-not}: for $A \in E^c$, we first use \eqref{eq:rosenthal-wp} and then \eqref{eq:nor-lp-hna} with $k=2$ and $q=1$, obtaining
\begin{equation}
\W_p\bra{  \bar{Z}_n(A),  \cM(A) +  \frac{\cN_{\cV(A)}}{\sqrt{n} } } \le \cost  \frac{\nor{Z(A) - \cM(A)}_{L^p }}{\sqrt{n}}  \le \cost \frac{ 1 + \nor{A} }{\sqrt{n}},
\end{equation}
with $\cost = \cost(p,S,T, h)<\infty$. Using
\begin{equation}
 \nor{A} \le \nor{A_0} + \nor{A-A_0} \le \cost \bra{1+ \nor{A-A_0}}
\end{equation}
inequality \eqref{eq:claim-prop-clt-square-E-not} follows, hence the thesis is proved.
\end{proof}

\subsection{Transportation of conditionally Gaussian variables}

We begin with the following simple observation: if $n$, $k \ge 1$ and $A$, $B$ are random variables taking values in $\R^{k \times k}$  and for $n \ge 1$, $\cN_{n\times k}$ denotes a standard Gaussian variable on $\R^{n \times k}$ that is independent of both $A$ and $B$, then
 \begin{equation}\label{eq:bound-abz}
  \W_p\bra{ (\Id_{n} \otimes A) \cN_{n \times k}, (\Id_{n} \otimes B) \cN_{n \times k}} \le \cost \sqrt{n} \W_p \bra{ A, B},
 \end{equation}
 where $\cost = \cost(p,k)< \infty$. Indeed, it is sufficient to consider any joint realization of $A$ and $B$ and a further independent variable $\cN_{n\times k}$ so that by \eqref{eq:moment-p-gaussian} we bound from above
 \begin{equation}
 \begin{split}
 W_p\bra{ (\Id_{n} \otimes A) \cN_{n \times k}, (\Id_{n} \otimes B) \cN_{n \times k} } & \le  \nor{ (\Id_{n} \otimes A) \cN_{n\times k} - (\Id_{n} \otimes B) \cN_{n \times k} }_{L^p}\\
 & =   \nor{ \bra{ \Id_{n} \otimes (A- B) } \cN_{n \times k} }_{L^p} \\
  & \le  \cost  \sqrt{ n }\nor{A-B}_{L^p},
  \end{split}
 \end{equation}
  where in the last identity we used the fact that the norm is multiplicative and $\nor{\Id_{n}} = \sqrt{n}$. Minimizing upon the joint realizations of $A$ and $B$ yields \eqref{eq:bound-abz}.

  Our aim in this section is to establish the following result which improves the dependence upon $\nor{A-B}_{L^p}$ in \eqref{eq:bound-abz}, when $B = B_0$ is constant and $A = B_0+G$, for an centered Gaussian variable $G$  with values in $\Sym^k$, independent of $\cN_{n\times k}$. This improvement however will come at the price of a worse dependence upon $n$.

\begin{proposition}\label{prop:gaussian-closer}
Let $k \ge 1$, $B_0 \in \Sym_+^k$ be invertible, let $G$
be a centered Gaussian random variable taking values in $\Sym^k$ with
\begin{equation}
 \nor{\Cov(G)} \le 1.
\end{equation}
For any $n \ge 1$, let $\cN_{n \times k}$ be a standard Gaussian random variable with values in $\R^{n \times k}$ and independent of $G$. Then,
\begin{equation}\label{eq:bound-gaussian-closer}
 \W_p \bra{\bra{ \Id_{n} \otimes (B_0+  G)} \cN_{n \times k}, \bra{ \Id_{n} \otimes B_0 } \cN_{n \times k}   } \le  \cost  n^2 \gamma_p^{nk}  \nor{\Cov(G)}
\end{equation}
where $\cost  = \cost(p,k, B_0)<\infty$ and $\gamma_p$ is defined in \eqref{eq:gamma-p}.
\end{proposition}

 Inequality \eqref{eq:bound-gaussian-closer} can be seen as a generalization of \cite[(5.14) in Proposition 5.2]{favaro2023quantitative} in the multivariate case and for the Wasserstein distance of any order $p\ge 1$. Let us remark that our proof is quite different from that in \cite{favaro2023quantitative} and does not rely upon Stein's method. Our argument  uses instead an upper bound for the Wasserstein distance in terms of suitable norm of the difference of the two probability measures. The following lemma, that is enough for our purposes,  is in fact a simplified version of a more general family of results (see e.g.\ \cite{peyre2018comparison, ambrosio2019pde, ledoux2019optimal,  goldman2021convergence, brown2020wasserstein}) that estimate from above the Wasserstein distance with a negative Sobolev norm. In our case we use Poincaré's inequality to further reduce the bound to an $L^p$ norm.

\begin{theorem}[transportation bound]\label{thm:ledoux}
 Let $X$ be a random variables with values in $\R^k$ and absolutely continuous law with respect to Lebesgue measure $\mathcal{L}^k$, with density
 \begin{equation}\label{eq:log-concave}
 \frac{d \mathbb{P}_X}{d \mathcal{L}^k} (x) = \exp\bra{-V(x)  - \frac{ \nor{x}^2}{2} },
 \end{equation}
 for some convex lower semicontinuous function $V: \R^k \to (-\infty, \infty]$. Then, for every $p \ge 1$ and any random variable $Y$ with values in $\R^k$ and absolutely continuous law with respect to that of $X$, it holds
 \begin{equation}
 \W_p(X, Y) = \cost \nor{ \frac{d\mathbb{P}_Y}{d\mathbb{P}_X}  (X)-  1}_{L^p},
\end{equation}
with $\cost  = \cost(p,k)<\infty$.
\end{theorem}
\begin{proof}
For brevity, we do not provide a detailed proof, which can be found e.g.\ in \cite[Lemma 3.4]{goldman2021convergence} for the case of a bounded rectangle endowed with Lebesgue measure and \cite[Theorem 2]{ledoux2019optimal} for an abstract case, which can be specialized on weighted Riemannian manifold with lower bounds on their curvature, in the sense of Bakry-\'Emery's $\Gamma_2$. Here, we only mention that the assumption \eqref{eq:log-concave} entails that the law of $X$ is $\log$-concave and actually with $\Gamma_2 \ge 1$, hence one can reproduce the general arguments (e.g.\ from \cite{ledoux2019optimal}) in this specialized setting, and apply well-known functional inequalities that are known to hold in this setting, such as Poincar\'e inequality and the $L^p$-boundedness of the Riesz transform (see also the recent work \cite{carbonaro2023boundedness}).
\end{proof}

Clearly, the above result applies in the case of $X = \cN_k$ being a standard Gaussian random variable. We will however apply it in the case of $X$ being a vector of independent $\chi_n$ laws: let us recall that the $\chi_n$ distribution (for $n \ge 1$) on the (half) real line is defined via the density
\begin{equation}
\chi_n(\rho) = \begin{cases} \frac{ \rho^{n-1} \exp\bra{-\rho^2/2} }{2^{n/2-1} \Gamma(n/2)} & \text{if $\rho \ge 0$,}\\
                 0 & \text{if $\rho <0$,}
               \end{cases}
\end{equation}
where $\Gamma$ denotes Euler's gamma function. Clearly, $\chi_n$ satisfies the assumptions of \Cref{thm:ledoux} with $V(\rho) = -(n-1) \log \rho$ if $\rho>0$, $V(\rho) = \infty$ otherwise, and therefore also any product measure on $\R^k$ consisting of $k\ge 1$ such measures.

Using \Cref{thm:ledoux}, we now establish the analogue of \eqref{prop:gaussian-closer} in the simpler case of $B_0 = \Id_k$ and $G$ being a uniform random variable on $\cur{A, -A}$ for some $A \in \Sym_+^k$ with small norm.

\begin{lemma}\label{lem:ledoux}
Let $k \ge 1$, $p \ge 1$ and $A \in \Sym_+^k$ with $\nor{A}_{op}\le 1/2$ and
\begin{equation}\label{eq:small-ledoux-norm-a}
  (1+\nor{A}_{op})^2  \le  \sqrt{p'},
\end{equation}
where $p' = p/(p-1) \in [1, \infty]$. Let $Z$ be a random variable uniformly distributed on $\cur{-1,1}$ and, for $n \ge 1$ let $\cN_{n \times k}$ be a standard Gaussian variable taking values in $\R^{n \times k}$, independent of $Z$. Then, it holds
\begin{equation}\label{eq:finite-n}
\W_p \bra{\Id_{n} \otimes \bra{ \Id_S + Z A }\cN_{n \times k}, \cN_{n \times k} } \le \cost \nor{A}^2 n^2 \gamma_p^{n k}
\end{equation}
where $\cost = \cost(p,k)< \infty$ depends upon $p$ and $k$ only and $\gamma_p$ is defined in \eqref{eq:gamma-p}.
\end{lemma}

\begin{proof}
We write for brevity $\cN = \cN_{n \times k}$. Using the spectral theorem, we decompose $A = U^\tau D U$ with $U$ orthogonal and $D$ diagonal, so that
\begin{equation}
 \Id_{n} \otimes \bra{ \Id_k+ Z A } \cN = (\Id_{n}\otimes U^\tau) \bra{  \Id_{n} \otimes \bra{ \Id_k + Z D} } (\Id_{n}\otimes U^\tau) \cN.
\end{equation}
Since  $\nor{\cdot}$ is unitarily invariant and
\begin{equation}
 (\Id_{n}\otimes U^\tau) \cN\stackrel{law}{=} \cN
\end{equation}
we obtain that
\begin{equation}
 \W_p\bra{ \Id_{n} \otimes \bra{ \Id_k+ Z A }\cN, \cN } = \W_p\bra{ \Id_{n} \otimes \bra{ \Id_k + Z D} \cN, \cN},
\end{equation}
hence we may reduce to the case of $A = D$ being diagonal. Writing $\cN = (\cN[:,i])_{i=1}^k$, where $\cN[:,i]$ are independent variables with values in $\R^n$ and standard Gaussian law, we can define,  for every $i =1, \ldots, k$,
\begin{equation}
 R_i := \nor{ \cN[:,i]}, \quad U_i := \frac{\cN[:,i]}{\nor{\cN[:,i]}},
\end{equation}
with the random variable $U_i$ taking values in $\R^{n}$ being uniformly distributed on the sphere $\mathbb{S}^{n-1}$ and independent of $R_i$. The variable $R = (R_i)_{i=1}^k$ takes values in $\R^k$ and consists of $k$ independent variables, each with $\chi_n$ law. Since $A$ is diagonal, we can also write
\begin{equation}
 \Id_{n} \otimes \bra{ \Id_k+ Z A } \cN= \bra{    \tilde{R}_i \tilde{U}_i}_{i=1}^k,
\end{equation}
where the variables $(\tilde{U}_i)_{i=1}^k$ have the same law as the variables $(U_i)_{i=1}^n$ and are further independent of the variables $\tilde{R} = (\tilde{R}_i)_{i=1}^k$. The variables $\tilde{R}_i$ however are not independent, since they have the same law as $\bra{ (1+a_i Z)R_i}_{i=1}^k$, where we write $a_i := A[i,i] \in [-1/2,1/2]$.
 Therefore, given any probability space where random variables $R = (R_i)_{i=1}^k$ and $\tilde{R} = \bra{\tilde{R}_i}_{i=1}^k$ are jointly defined (each with the law just described), up to enlarging the space we can consider a further independent collection of $k$ independent variables $\bra{ U_i}_{i=1}^k$ uniformly distributed on the sphere $\mathbb{S}^{n-1}$ to define, for $i=1, \ldots, k$,
\begin{equation}
 N_i := R_i U_i, \quad \text{and}\quad \tilde{N}_i :=\tilde{R}_i U_i,
\end{equation}
so that the resulting vectors $N = (N_i)_{i=1}^k  \stackrel{law}{=} \cN$ has standard Gaussian distribution and $\tilde{N} = (\tilde{N}_i)_{i=1}^k$ has the same law as $(\Id_n \otimes \bra{ \Id_k+ Z A } \cN_{n\times k})$. Using this construction, by \eqref{eq:w-2-trivial-bound} we bound from above
\begin{equation}\begin{split}
\W_p\bra{ \Id_{n} \otimes \bra{ \Id_k+ Z A } \cN,  \cN }  & \le  \nor{ \tilde N-N}_{L^p} \le \nor{ \bra{ \sum_{i=1}^k \nor{ (\tilde{R}_i - R_i) U_i}^2 }^{1/2}}_{L^p}\\
& = \nor{  \bra{ \sum_{i=1}^k \abs{\tilde{R}_i - R_i}^2  }^{1/2} }_{L^p}  =  \nor{ \tilde R - R}_{L^p}
\end{split}
\end{equation}
Minimizing upon the joint realization of $R$ and $\tilde{R}$ yields the inequality
\begin{equation}
 \W_p\bra{ \Id_{n} \otimes \bra{ \Id_k+ Z A } \cN,  \cN }  \le \W_p\bra{ \tilde{R}, R},
\end{equation}
which can be rewritten as
\begin{equation}
  \W_p\bra{ \Id_{n} \otimes \bra{ \Id_k+ Z A } \cN,  \cN }  \le \W_p\bra{ (1+Z A) R, R},
\end{equation}
where $R = (R_i)_{i=1}^k$ are independent variables with law $\chi_n$ and $Z$ is independent of $R$ and uniform on $\cur{-1,1}$. Since $\nor{A}_{op} \le 1/2 <1$, the variable $(1+Z A) R$ has absolutely continuous law with respect to that of $R$ and is given  by following the expression
 \begin{equation}\label{eq:f-density}
  f((r_i)_{i=1}^k, (a_i)_{i=1}^k ) = \frac 1 2  \prod_{i=1}^k \frac{\exp\bra{\frac{r_i^2}{2} \bra{1-\frac{1}{(1+a_i)^2} }}}{(1+a_i)^n}  +  \frac 1 2\prod_{i=1}^k \frac{\exp\bra{\frac{r_i^2}{2} \bra{1-\frac{1}{(1-a_i)^2} } }}{(1-a_i)^n}.
 \end{equation}
Therefore, by \Cref{lem:ledoux}, we bound from above
\begin{equation}\label{eq:ledoux-1sa}
 \W_p\bra{ (1+Z A) R, R} \le \cost \nor{ f(R, (a_i)_{i=1}^k )- 1}_{L^p}.
\end{equation}
For any given $r = (r_i)_{i=1}^k$, the function $a =(a_i)_{i=1}^k \mapsto f(r,a)$ is smooth on $[-1/2,1/2]^k$, it holds $f(r,0) = 1$  and is even: $f(r,a) = f(r,-a)$.  Therefore, if we consider the smooth and even function $h(r, t) = f(r, ta)$ for $t \in [-1,1]$ it holds $h(r, 0) = 1$ and $\partial_t h(r,0) = 0$. Hence, by Taylor's expansion, we have
\begin{equation}\label{eq:taylor-h}
 f(r,  a ) - 1 = h(r, 1) - h(r, 0) =  \frac 1 2 \int_0^1 \partial^2_t h(r, s) (1-s) ds.
\end{equation}
Since
\begin{equation}
 \partial^2_t h(r,s)  = \sum_{i,j=1}^k \frac{\partial^2 f}{\partial a_i \partial a_j} (r, sa) a_i a_j
\end{equation}
we conclude that
\begin{equation}
 f(r, a ) - 1  = \int_0^1 \sum_{i,j=1}^k \frac{\partial^2 f}{\partial a_i \partial a_j} (r, sa)  a_i a_j (1-s) ds
\end{equation}
Using this identity in the right hand side of \eqref{eq:ledoux-1sa} yields
\begin{equation}
\begin{split}
\nor{ f(R, a )- 1}_{L^p} & \le \int_0^1 \sum_{i,j=1}^k \abs{a_i a_j} \nor{ \frac{\partial^2 f}{\partial a_i \partial a_j} (R, sa)}_{L^p} ds\\
 & \le \bra{\sum_{i=1}^k |a_i|}^2 \sup_{s \in [0,1], i,j=1, \ldots, k} \nor{  \frac{\partial^2 f}{\partial a_i \partial a_j} (R, sa)}_{L^p}.
 \end{split}
\end{equation}
Thus, to obtain the thesis we only need to prove that
\begin{equation}
  \nor{  \frac{\partial^2 f}{\partial a_i \partial a_j} (R, (\alpha_\ell)_{\ell=1}^k) }_{L^p} \le \cost n^2 \gamma_p^{nk},
\end{equation}
if $\max_{\ell=1, \ldots, k}{|\alpha_\ell|} \le \nor{A}$, hence by \eqref{eq:small-ledoux-norm-a} it holds $(1+|\alpha_\ell|)^2 \le \sqrt{p'}$ for every $\ell=1,\ldots, k$. To this aim, we go back to the definition of $f$ in \eqref{eq:f-density} and rewrite it as
\begin{equation}
 f(r, (\alpha_{\ell})_{\ell=1}^k ) = \frac 1 2 \prod_{\ell=1}^k g(r_\ell, \alpha_\ell) + \frac 1 2 \prod_{\ell=1}^k g(r_\ell, -\alpha_\ell),
\end{equation}
where
\begin{equation}
 g(\rho,\alpha) = \frac{1}{(1+\alpha)^n} \exp\bra{ \frac{\rho^2}{2} (1- (1+\alpha)^{-2})}.
\end{equation}
By Minkowski inequality for the $L^p$ norm, we find
\begin{equation}
 \nor{\frac{\partial^2 f}{\partial a_i \partial a_j} (R, (\alpha_\ell)_{\ell=1}^{k}) }_{L^p} \le \frac 1 2  \nor{\frac{\partial^2 \prod_{\ell =1 }^k g  (R_\ell, \alpha_\ell)}{\partial a_i \partial a_j} }_{L^p} + \frac 1 2\nor{\frac{\partial^2 \prod_{\ell=1}^k g(R_\ell, -\alpha_\ell)}{\partial a_i \partial a_j} }_{L^p}
\end{equation}
and focus for simplicity on the first term (the second one is similar).

We notice the identities
\begin{equation}\begin{split}
\frac{ \partial g }{\partial \alpha} (\rho, \alpha)& = \sqa{ - n (1+\alpha)^{-1}  +  \rho^2(1+\alpha)^{-3}}  g(\rho, \alpha)\\
& := m_1(\rho, \alpha) g(\rho, \alpha)
\end{split}
\end{equation}
and
\begin{equation}
\begin{split}
\frac{ \partial^2 g }{\partial \alpha^2} (\rho, \alpha) & = \sqa{n (1+\alpha)^{-2} -3 \rho^2(1+\alpha)^{-4} + m_1(\rho, \alpha)} g(\rho, \alpha)\\
& := m_2(\rho, \alpha) g(\rho, \alpha).
\end{split}
\end{equation}
Using these identities we see that, for $i \neq j$,
\begin{equation}
 \frac{\partial^2 \prod_{\ell=1}^k g(r_\ell, \alpha_\ell)}{\partial a_i \partial a_j} =  m_1(r_i, \alpha_i)m_1(r_j, \alpha_j) \prod_{\ell=1}^k g(r_\ell, \alpha_\ell),
\end{equation}
while for  $i= j$,
\begin{equation}
 \frac{\partial^2 \prod_{\ell=1}^k g(r_\ell, \alpha_\ell)}{\partial a_i^2} =  m_2(r_i, \alpha_i) \prod_{\ell=1}^k g(r_\ell, \alpha_\ell).
\end{equation}
In view of the independence of the variables $(R_i)_{i=1}^k$ and the product form of the functions above, it is sufficient to  bound from above the quantities
\begin{equation}
 \nor{ g(R_1, \alpha)}_{L^p}, \quad \nor{ m_1(R_1, \alpha) g(R_1, \alpha)}_{L^p} \quad \text{and} \quad  \nor{ m_2(R_1, \alpha) g(R_1, \alpha)}_{L^p}.
\end{equation}
Consider the case $p=1$ first. Then, $g(\rho, \alpha)$ is the density of the random variable $(1+\alpha)R_1$ with respect to the law of $R_1$ (that is $\chi_n$), hence for every $\varphi: (0, \infty) \to [0, \infty)$, we have the change of variables
\begin{equation}
 \nor{ \varphi(R_1) g(R_1, \alpha)}_{L^1} = \int_0^\infty \varphi(\rho) g(\rho, \alpha) \chi_n(\rho) d \rho = \EE\sqa{ \varphi((1+\alpha)R_1)}.
 \end{equation}
Using $\varphi(\rho) = 1$ yields $ \nor{ g(R_1, \alpha)}_{L^1} =1$, while $\varphi(\rho) = m_1(\rho, \alpha)$ gives
\begin{equation}
  \nor{ m_1(R_1, \alpha) g(R_1, \alpha)}_{L^1} \le \cost n,
\end{equation}
where $\cost < \infty$, using the fact that we are assuming $|\alpha| \le 1/2$ and that the second moment of a $\chi_n$ law equals $n$. Similarly,
\begin{equation}
 \nor{ m_2(R_1, \alpha) g(R_1, \alpha)}_{L^1} \le \cost n.
\end{equation}
These inequalities lead easily to the bound
\begin{equation}
  \nor{\frac{\partial^2 \prod_{\ell =1 }^k g  (R_\ell, \alpha_\ell)}{\partial a_i \partial a_j} }_{L^1} \le \cost n^2,
\end{equation}
and eventually the thesis.

For the case, $p>1$, we  simply observe that
\begin{equation}\begin{split}
  g(\rho, \alpha)^p \chi_n(\rho)  & = (1+\alpha)^{-n} \exp\bra{ - \frac{\rho^2}{2} (1 +p((1+\alpha)^{-2}-1))} \frac{ \rho^{n-1} }{2^{n/2-1}\Gamma(n/2)} \\
  & = \bra{ \frac{ \beta}{1+\alpha}}^n \exp\bra{ - \frac{(\rho/\beta)^2}{2}  } \frac{ (\rho/\beta)^{n-1} }{2^{n/2-1}\Gamma(n/2) \beta},
  \end{split}
\end{equation}
where the parameter $\beta$ is given by the expression
\begin{equation}
 \beta := \frac{ 1}{\sqrt{1+p((1+\alpha)^{-2}-1)}} = \frac{1+\alpha}{\sqrt{p} \sqrt{1 -(1+\alpha)^2/p' }}
\end{equation}
Notice that $\beta$ is well defined, provided that $(1+\alpha)^2 < p'$, but in view of \eqref{eq:small-ledoux-norm-a} it actually holds
\begin{equation}
 (1+\alpha)^2 < \sqrt{ p'}
\end{equation}
hence, using that for every $x \in (0,1)$ it holds $\sqrt{1-x} \ge 1-|x|$, we also obtain the inequality
 \begin{equation}
  \frac{\beta}{1+\alpha} \le \frac{ 1}{ \sqrt{p} \bra{ 1- 1/\sqrt{p'} }}= \frac{1}{\sqrt{p} - \sqrt{p-1}} = \gamma_p^p.
 \end{equation}
Thus,  we find that for every function $\varphi: (0, \infty) \to [0, \infty)$,
\begin{equation}\begin{split}
 \nor{ \varphi(R_1)  g(R_1, \alpha)}_{L^p}^p & = \int_0^\infty \varphi(\rho)^p g(\rho, \alpha)^p \chi_n(\rho) d \rho \\
 & \le \gamma_p^{pn}  \EE\sqa{ \varphi(\beta R_1)^p},
 \end{split}
\end{equation}
Therefore, we have trivially
\begin{equation}
 \nor{ g(R_1, \alpha) }_{L^p}\le \gamma_p^{n}
\end{equation}
and similarly we find
\begin{eqnarray}
 \nor{ m_1(R_1, \alpha) g(R_1, \alpha)}_{L^p} & \le \cost n \gamma_p^{n},\\
  \nor{ m_2(R_1, \alpha) g(R_1, \alpha)}_{L^p} & \le \cost n \gamma_p^n.
\end{eqnarray}
Combining these bounds, we obtain
\begin{equation}
  \nor{\frac{\partial^2 \prod_{\ell =1 }^k g  (R_\ell, \alpha_\ell)}{\partial a_i \partial a_j} }_{L^p} \le \cost n^2 \gamma_p^{nk},
\end{equation}
and eventually the thesis also in this case.
\end{proof}

We are now in a position to prove \Cref{prop:gaussian-closer}.

\begin{proof}[Proof of \Cref{prop:gaussian-closer}]
For brevity, we  write $\cN = \cN_{n\times k}$ in what follows. Moreover, writing
\begin{equation}
 \Id_{n} \otimes (B_0+  G)  = (\Id_{n} \otimes B_0) \bra{ \Id_{n} \otimes (\Id_k+ B_0^{-1}G)},
\end{equation}
we bound from above
\begin{equation}
\begin{split}
& \W_p \bra{\bra{ \Id_{n} \otimes (B_0+  G)} \cN, \bra{ \Id_{n} \otimes B_0 } \cN   } \le\\
& \quad \le \nor{(\Id_{n} \otimes B_0) }_{op}\W_p \bra{\bra{ \Id_{n} \otimes (\Id_k+  B_0^{-1} G)}  \cN,   \cN   }\\
& \quad =\nor{B_0}_{op} \W_p \bra{\bra{ \Id_{n} \otimes (\Id_k+  B_0^{-1} G)}  \cN,   \cN   }.
\end{split}
\end{equation}
Hence, up to replacing $G$ with $B_0^{-1} G$ we can assume that $B_0 = \Id_k$ and write for simplicity $G$ instead of $B_0^{-1}G$. Notice however that now $G$ is not symmetric (but still Gaussian).  Still, we have the identity in law
\begin{equation}
 (\Id_{n} \otimes (\Id_k+ G) ) \cN \stackrel{law}{=} (\Id_{n} \otimes  | \Id_k + G | )  \cN
\end{equation}
where $|\Id_k + G| = \sqrt{(\Id_k + G)(\Id_k + G^\tau)}$, and using the triangle inequality,
\begin{equation}\label{eq:ledoux-first-second-terms}\begin{split}
 &  \W_p \bra{ (\Id_{n} \otimes  |\Id_k+G|) \cN,   \cN   } \le \\
  & \le \W_p \bra{  (\Id_{n} \otimes  |\Id_k +G|) \cN,   \Id_{n} \otimes  \bra{ \Id_k+ G^{\sym} }  \cN  }\\
 & \quad  + \W_p \bra{ \bra{\Id_{n} \otimes \bra{ \Id_k+ G^{\sym}} } \cN,  \cN }
 \end{split}
\end{equation}
where we let $G^{\sym} := (G +G^{\tau})/2$. We focus separately on the two terms: for the first term, by \eqref{eq:bound-abz}, we bound from above
\begin{equation}\begin{split}
  & \W_p \bra{ \bra{ \Id_{n} \otimes |\Id_k + G| }  \cN,   \bra{\Id_{n} \otimes  \bra{ \Id_k+ G^{\sym} }}  \cN  } \le\\
  & \le \cost \sqrt{n} \W_p \bra{  |\Id_k + G| ,   \Id_k+ G^{\sym}  }.
  \end{split}
\end{equation}
We apply \Cref{prop:delta-method-sqrt} with $A_0 := \Id_{k}$ and
\begin{equation}
 A_1:= |\Id_{k}+ G|^2 = \Id_k + 2 G^{\sym}+ |G|^2
\end{equation}
so that $A_1 - A_0 = 2G^{\sym}+|G|^2$, and $D_{\sqrt{\cdot}}(\Id_k)B = \frac 1 2 B$, this leads to the inequality
\begin{equation}
\W_p \bra{  |\Id_k + G| , \Id_k + G^{\sym} + |G|^2 } \le \cost \nor{ G + G^\tau + |G|^2 }^2_{L^{2p}} \le \cost \nor{\Cov(G)},
                \end{equation}
having used that $\nor{\Cov(G)} \le 1$. Using \eqref{eq:w-2-trivial-bound}, we also trivially bound from above
\begin{equation}
 \W_p\bra{ \Id_k + G^{\sym} + |G|^2 , \Id_k+  G^{\sym} } \le \nor{|G|^2}_{L^p}  \le \cost \nor{\Cov(G)},
\end{equation}
with $\cost = \cost(p,k)< \infty$. Hence, again by the triangle inequality, we conclude that
\begin{equation}\begin{split}
 & \W_p \bra{  |\Id_k + G| ,   \Id_k+ G^{\sym}  } \le \\
 & \le \W_p \bra{  |\Id_k + G| , \Id_k +G^{\sym} + |G|^2  } \\
 & \quad + \W_p \bra{ \Id_k + G^{\sym} + |G|^2,  \Id_k+ G^{\sym}} \\
  & \le \cost \nor{ \Cov(G)}.
  \end{split}
\end{equation}
with $\cost = \cost(p,k)<\infty$. We deduce that the contribution of the first term in \eqref{eq:ledoux-first-second-terms} is
\begin{equation}\label{eq:bound-first-term-lemma-ledoux}
\begin{split}
 & \W_p \bra{  \bra{\Id_{n} \otimes |\Id_k+  G |} \cN,  \bra{ \Id_{n} \otimes  \bra{ \Id_k+ G^{\sym} }} \cN  } \le \\
 & \le \cost \sqrt{n} \nor{ \Cov(G)}.
 \end{split}
 \end{equation}
For the second term in \eqref{eq:ledoux-first-second-terms}, we claim that it holds
\begin{equation} \label{eq:bound-second-term-lemma-ledoux}
 \W_p \bra{\Id_{n} \otimes \bra{ \Id_k+ G^{\sym}} \cN, \cN } \le \cost n^2 \gamma_p^{nk} \nor{\Cov (G)},
\end{equation}
which  combined with \eqref{eq:bound-first-term-lemma-ledoux} yields the thesis. To simplify the notation, we write now $G$ instead of $G^{\sym}$ (and will use that $G$ is symmetric). We introduce a further independent random variable $Z$, with uniform law in $\cur{-1, 1}$, so that by the assumption on the law of $G$ being equal to that of $-G$, it follows  that also
\begin{equation}
  Z G \stackrel{law}{=} G.
\end{equation}
Thus,
\begin{equation}
  \W_p \bra{ \bra{\Id_{n} \otimes \bra{ \Id_k+ G}} \cN, \cN } =  \W_p \bra{ \bra{ \Id_{n} \otimes \bra{ \Id_k+ ZG}} \cN, \cN }.
\end{equation}
Next, we invoke the convexity inequality \eqref{eq:convexity}, so that
\begin{equation}
\W_p \bra{ \bra{\Id_{n} \otimes \bra{ \Id_k+ Z G}} \cN, \cN } \le \bra{ \int_{\R^{k\times k}} \W_p \bra{\Id_{n} \otimes \bra{ \Id_k+ Z A} \cN, \cN }^p d \mathbb{P}_{G}(A)}^{1/p}
\end{equation}
where $\mu$ denotes the law of $G (=G^{sym})$. We introduce the set
\begin{equation}
 E = \cur{A\,  : \,  \nor{A}_{op} \le 1/2 \text{ and }  (1+\nor{A}_{op})^2 \le \sqrt{p'}}
\end{equation}
on which we can apply  \Cref{lem:ledoux}, yielding
\begin{equation}
  \W_p \bra{\Id_{n} \otimes \bra{ \Id_k+ Z A} \cN, \cN } \le \cost n^2 \gamma_p^{nk} \nor{A}^2
\end{equation}
If instead $A \in E^c$, we use \eqref{eq:bound-abz} to bound from above
\begin{equation}
 \W_p \bra{ \bra{\Id_{n} \otimes \bra{ \Id_k+ Z A} }\cN, \cN }  \le \cost \sqrt{n} \nor{A}.
\end{equation}
Splitting the integration over $E$ and $E^c$, we obtain therefore
\begin{equation}\begin{split}
& \bra{ \int_{\R^{k\times k}} \W_p\bra{\Id_{n} \otimes \bra{ \Id_k+ Z A} \cN, \cN }^p d \mathbb{P}_{G}(A)}^{1/p} \le \\
& \le \bra{ \int_{E} \W_p \bra{\Id_{n} \otimes \bra{ \Id_k+ Z A} \cN, \cN }^p d \mathbb{P}_{G}(A)}^{1/p} \\
& \quad + \bra{ \int_{E^c} \W_p \bra{\Id_{n} \otimes \bra{ \Id_k+ Z A} \cN, \cN }^p d \mathbb{P}_{G}(A)}^{1/p}\\
& \le \cost \bra{ n^2 \gamma_p^{nk} \EE\sqa{ \nor{G}^{2p}}^{1/p} + \sqrt{n} \EE\sqa{ I_{\cur{ G \in E^c}} \nor{G}^{p}}^{1/p}}
\end{split}
\end{equation}
To conclude, we notice that on $E^c$ it holds $\nor{A}_{op} > 1/2$ or $\nor{A}_{op} \ge (p')^{1/4} -1$, hence $\nor{A}_{op} \ge \eps$ for some $\eps = \eps(p)< \infty$. Therefore, by Markov's inequality, we have
\begin{equation}\label{eq:markov-ec-ledoux}
 \mathbb{P} (G \in E^c) \le \frac{\nor{G}_{L^{2p}}^{2p}}{\eps^{2p}} \le \cost \nor{\Cov(G)}^{p},
\end{equation}
with $\cost = \cost(p,q, k) <\infty$, hence by Cauchy-Schwarz inequality we find
\begin{equation}
 \EE\sqa{ I_{ \cur{ G \in E^c}} \nor{G}^{p}}^{1/p} \le  \mathbb{P} (G \in E^c)^{1/(2p)} \nor{G}_{L^{2p}} \le \cost \nor{ \Cov(G)},
\end{equation}
and the proof is complete.
\end{proof}

\section{Main Result}\label{sec:proof}

Throughout this section, we consider the generalized fully connected architecture introduced in
\Cref{sec:nn} for a given  input set $\X = \cur{x_i}_{i=1}^k \subseteq \R^{d_0}$ and auxiliary dimensions $\bm{S}$, $\bm{a}$ and activation functions $\bm{\sigma}$. We then let the weights ${\bf W}$  to be independent random variables with standard Gaussian laws rescaled by the square root of the respective input layer widths, i.e., for $\ell=1, \ldots, L$,
\begin{equation}\label{eq:law-weights}
 W^{(\ell)} := \frac 1 {\sqrt{n_{\ell-1}}} \cN_{n_{\ell} \times (n_{\ell-1} \times a_{\ell})}
\end{equation}
where we recall that we also set for convenience $n_0 :=1$ and $a_0:= d_0 \times k$.

\subsection{Quantitative approximation in the possibly degenerate case}

As a warm-up and for a useful comparison with our improved rates (\Cref{thm:main-induction} below), we first establish the following result, borrowing most of the arguments from \cite{basteri2022quantitative}.

\begin{theorem}\label{thm:main-induction-basteri}
With the notation above, assume that all the activation functions $\bm{\sigma} = (\sigma^{(\ell)})_{\ell = 1}^L$ are Lipschitz continuous. Then, for every $p \ge 1$, there exists a constant
\begin{equation}
 \cost = \cost(p, \X, \bm{S}, \bm{a}, \bm{\sigma}, \bm{K})< \infty
\end{equation}
not depending on the hidden and output layer widths $\bm{n} = (n_i)_{i=1}^L$, such that, for every $\ell \in \cur{0,1,\ldots, L}$, it holds
\begin{equation}\label{eq:main-induction-baster} \W_p\bra{ f^{(\ell)}, G^{(\ell)}_{n_\ell}} \le \cost  \sqrt{n_\ell} \sum_{\ell =1}^{\ell-1}\frac 1 {\sqrt{ n_i}}.\end{equation}
and for every Lipschitz continuous function $h: \R^{S_\ell} \to \R^{T}$ (with $T$ any finite set) it holds
\begin{equation}
 \label{eq:induction-kernel-basteri}
\W_p\bra{  \frac 1 {n_\ell} \sum_{i=1}^{n_\ell} \bra{h(f^{(\ell)}[i,:])}^{\otimes 2}, \EE\sqa{ \bra{ h(G^{\ell}) }^{\otimes 2} } }   \le  \cost \nor{h}^2_{\Lip}   \sum_{i=1}^{\ell} \frac{1}{\sqrt{n_i}}.
\end{equation}
\end{theorem}

\begin{proof}
We argue by induction over $\ell \in \cur{0,1,\ldots, L}$ to establish the validity of \eqref{eq:main-induction-baster} and \eqref{eq:induction-kernel-basteri}. Notice also that we can assume $p \ge 2$ without loss of generality. The case $\ell=0$ is trivially true, since by construction we have $f^{(0)} = G^{(0)}_{n_0}$ is constant (and $n_0:=1$).

Next, we assume that the thesis holds for $\ell-1$ and establish its validity for the case $\ell$. We prove first \eqref{eq:main-induction-baster}. By the triangle inequality, we split
\begin{equation}\label{eq:basteri-first-second-term}\begin{split}
  \W_p\bra{ f^{(\ell)}, G^{(\ell)}_{n_\ell}} & \le \W_p\bra{f^{(\ell)},   (W^{(\ell)} \otimes \Id_{S_\ell})  \sigma^{(\ell)} ( G^{(\ell-1)}_{n_{\ell-1}} ) } \\
  & \quad + \W_p\bra{(W^{(\ell)} \otimes \Id_{S_\ell})  \sigma^{(\ell)} ( G^{(\ell-1)}_{n_{\ell-1}} ), G^{(\ell)}_{n_\ell}}
  \end{split}
\end{equation}
where $W^{(\ell)}$ is independent of $G^{(\ell-1)}_{n_{\ell-1}}$, and we bound separately the two terms. For the first term, we recall that
\begin{equation}
 f^{(\ell)}  = (W^{(\ell)}\otimes \Id_{S_\ell}) \sigma^{(\ell)}(f^{(\ell-1)})
\end{equation}
hence by \eqref{eq:w-2-trivial-bound} we find
\begin{equation}\begin{split}
& \W_p\bra{f^{(\ell)},   (W^{(\ell)} \otimes \Id_{S_\ell})  \sigma^{(\ell)} ( G^{(\ell-1)}_{n_{\ell-1}} ) }  \le\\
& \le \nor{ \bra{ (W^{(\ell)} \otimes \Id_{S_\ell} } \bra{   \sigma^{(\ell)}(f^{(\ell-1)}) -  \sigma^{(\ell)} ( G^{(\ell-1)}_{n_{\ell-1}} ) } }_{L^p}\\
& \le \EE\sqa{ \EE\sqa{ \nor{ (W^{(\ell)} \otimes \Id_{S_\ell} ) v }^p \bigg| v = \sigma^{(\ell)}(f^{(\ell-1)}) -  \sigma^{(\ell)} ( G^{(\ell-1)}_{n_{\ell-1}} } }^{1/p}\\
&  \le \cost \sqrt{n_\ell} \nor{ \sigma^{(\ell)}(f^{(\ell-1)}) -  \sigma^{(\ell)} ( G^{(\ell-1)}_{n_{\ell-1}}) }_{L^p}
\end{split}
\end{equation}
where $\cost = \cost(p, a_\ell)$, having conditioned upon the value of $\sigma^{(\ell)}( f^{(\ell-1)})-\sigma^{(\ell)}( G^{(\ell-1)}_{n_{\ell-1}})$. Since $\sigma^{(\ell)}$ acts componentwise and is assumed to be Lispschitz continuous, it follows that
\begin{equation}
 \nor{ \sigma^{(\ell)}(f^{(\ell-1)}) -  \sigma^{(\ell)} ( G^{(\ell-1)}_{n_{\ell-1}}) }_{L^p} \le \Lip(\sigma^{(\ell)}) \nor{ f^{(\ell-1)} -  G^{(\ell-1)}_{n_{\ell-1}}}_{L^p}.
\end{equation}
Minimizing upon the join realizations of $f^{(\ell-1)}$ and $G^{(\ell-1)}$ yields by the inductive assumption
\begin{equation}
 \W_p\bra{f^{(\ell)},   (W^{(\ell)} \otimes \Id_{S_\ell})  \sigma^{(\ell)} ( G^{(\ell-1)}_{n_{\ell-1}} ) } \le \cost \sqrt{n_{\ell}} \sum_{i=1}^{\ell-1}\frac 1 {\sqrt{n_i}},
 \end{equation}
 hence  to obtain \eqref{eq:main-induction-baster} we are only left with estimating the second term in \eqref{eq:basteri-first-second-term}.  To this aim, we use \eqref{eq:identity-in-law} with $A := \sigma^{(\ell)} ( G^{(\ell-1)}_{n_{\ell-1}} )$ $S:=S_\ell$, $T_2:= n_\ell$, $T_1:= n_{\ell-1} \times a_\ell$, so that
  \begin{equation}
   (W^{(\ell)} \otimes \Id_{S_\ell})  \sigma^{(\ell)} ( G^{(\ell-1)}_{n_{\ell-1}}) \stackrel{law}{=}  \cN_{ \Id_{n_\ell} \otimes \Sigma }
 \end{equation}
  where we write for brevity
 \begin{equation}
   \Sigma:= \frac 1 {n_{\ell-1}} \tr_{n_{\ell-1}\times a_\ell}(\sigma^{(\ell)} ( G^{(\ell-1)}_{n_{\ell-1}})^{\otimes 2}).
  \end{equation}
  By construction, we also have
  \begin{equation}
    G^{(\ell)}_{n_\ell} \stackrel{law}{=}\cN_{ \Id_{n_\ell} \otimes K^{(\ell)} }.
  \end{equation}
  Therefore can use \eqref{eq:convexity} with $U = \Sigma$ and then \eqref{eq:gaussian} to obtain
 \begin{equation}\begin{split}
  \W_p\bra{  (W^{(\ell)} \otimes \Id_{S_\ell})  \sigma^{(\ell)} ( G^{(\ell-1)}_{n_{\ell-1}}), G^{(\ell)}_{n_\ell}} &  \le \cost\nor{ \sqrt{\Id_{n_\ell} \otimes \Sigma} - \sqrt{ \Id_{n_\ell} \otimes K^{(\ell)}}}_{L^p}\\
  & \le \cost\sqrt{n_\ell} \nor{ \sqrt{\Sigma} - \sqrt{  K^{(\ell)}}}_{L^p}.
  \end{split}
 \end{equation}
We rewrite
\begin{equation}\label{eq:sigma-sum-iid}
 \Sigma =   \frac 1 {n_{\ell-1}} \sum_{i=1}^{n_{\ell-1}} \tr_{a_\ell}(\sigma^{(\ell)} ( G^{(\ell-1)}_{n_{\ell-1}}[i,:])^{\otimes 2}),
\end{equation}
and since $G^{(\ell-1)}_{n_{\ell-1}}[i,:] \stackrel{law}{=} G^{(\ell-1)}$ for every $i$, recalling the definition of $K^{(\ell)}$ we have
\begin{equation}
 \EE\sqa{ \tr_{a_\ell} \bra{ \sigma^{(\ell)} ( G^{(\ell-1)}_{n_{\ell-1}}[i,:])^{\otimes 2} } } =  K^{(\ell)}
\end{equation}
for every $i=1,\ldots, n_{\ell-1}$. By \eqref{eq:Xtwo-ac} it follows that
\begin{equation}
 \tr_{a_\ell} \bra{ \sigma^{(\ell)} ( G^{(\ell-1)}_{n_{\ell-1}}[i,:])^{\otimes 2} } \ll  K^{(\ell)},
\end{equation}
hence  $\Sigma \ll K^{(\ell)}$.  Therefore, we are in a position to apply \eqref{eq:ando}, yielding
\begin{equation}
 \nor{ \sqrt{\Sigma} - \sqrt{ K^{(\ell)}}}_{L^p} \le \cost \nor{\Sigma -  K^{(\ell)}}_{L^p},
\end{equation}
with $\cost = \cost(p,K^{(\ell)})< \infty$ (notice that if $K^{(\ell)} = 0$ the inequality trivializes). To conclude, it is sufficient to use Rosenthal's inequality \eqref{eq:rosenthal} with $n:= n_{\ell-1}$ and $d  =|S_\ell|$ and the  independent (centered) variables
\begin{equation}
 Z_i  = \tr_{a_\ell}(\sigma^{(\ell)} ( G^{(\ell-1)}_{n_{\ell-1}}[i,:])^{\otimes 2}) - K^{(\ell)}
\end{equation}
for $i=1, \ldots, n_\ell$. We obtain
\begin{equation}
\nor{\Sigma -  K^{(\ell)}}_{L^p} \le  \cost \frac{\nor{ Z_i}_{L^p} }{\sqrt{n_{\ell-1}}} \le \cost \frac{1}{\sqrt{n_{\ell-1}}},
\end{equation}
because using that $\sigma^{(\ell)}$ is Lipschitz (hence it has polynomial growth of order $q=1$) it follows easily that
\begin{equation}
 \nor{Z_i}_{L^p} \le \cost  = \cost(p,a_\ell, \sigma^{(\ell)}, K^{(\ell-1)}, K^{(\ell)} )< \infty.
\end{equation}
Hence, also the second term in \eqref{eq:basteri-first-second-term} is properly bounded from above and \eqref{eq:main-induction-baster} follows in the case $\ell$.

To establish \eqref{eq:induction-kernel-basteri}, without loss of generality we may assume $\nor{h}_{\Lip} \le 1$. We claim that the thesis is a consequence of the general bound \eqref{eq:kernel-basteri-tool} with $X:= f^{(\ell)}$ and $Y := G^{(\ell)}_{n_\ell}$  and using \eqref{eq:main-induction-baster} in the case $\ell$ (and $2p$ instead of $p$). Indeed, we find immediately that
\begin{equation}\begin{split}
 W_p\bra{   \sum_{i=1}^{n_\ell} (h(f^{(\ell)}[i,:]))^{\otimes 2},  \sum_{i=1}^{n_\ell} (h(G_{n_\ell}^{(\ell)}[i,:]))^{\otimes 2}} &  \le \cost \bra{ 1 + \sqrt{n_\ell} \sum_{i=1}^{\ell-1} \frac 1 {\sqrt{n_i}}} \sqrt{n_\ell} \sum_{i=1}^{\ell-1} \frac 1 {\sqrt{n_i}}\\
 & \le \cost n _{\ell} \sum_{i=1}^{\ell-1} \frac 1 {\sqrt{n_i}}.
 \end{split}
\end{equation}
Next, by \eqref{eq:w-2-trivial-bound} and \eqref{eq:rosenthal} with $n := n_\ell$ applied to the random variables
\begin{equation}
 Z_i := (h(G_{n_\ell}^{(\ell)}[i,:]))^{\otimes 2}
\end{equation}
which are all independent with the same law as $Z:=(h(G_{n_\ell}))^{\otimes 2}$, we deduce that
\begin{equation}
 \W_p\bra{\frac 1 {n_\ell} \sum_{i=1}^{n_\ell} (h(G_{n_\ell}^{(\ell)}[i,:]))^{\otimes 2}, \EE\sqa{h(G^{(\ell)})^{\otimes 2} } } \le \cost \frac{\nor{ Z}_{L^p}}{\sqrt{n_\ell}}  \le  \cost \cdot \frac{1}{\sqrt{n_\ell}},
\end{equation}
having finally used \eqref{eq:nor-lp-hna} with $k=2$, $q=1$ and $A:= K^{(\ell-1)}$. By the triangle inequality for $\W_p$ the validity of \eqref{eq:induction-kernel-basteri} follows.
\end{proof}

 \subsection{Quantitative approximation in the non-degenerate case}

 We are now in a position to state and prove the main result of this work, which provides improved rates for the Gaussian approximation of a random deep neural network assuming non-degeneracy of the infinite-width covariances on a (finite) input set.

\begin{theorem}\label{thm:main-induction}
With the notation above (in particular those of \Cref{sec:nn}, \eqref{eq:law-weights} and \eqref{eq:gamma-p}),  assume that all the activation functions $\bm{\sigma} = (\sigma^{(\ell)})_{\ell = 1}^L$ are Lipschitz continuous and the infinite-width covariances $\bm{K}$ are non degenerate on the input set $\X$.

Then, for every $p \ge 1$, there exists a constant
\begin{equation}
 \cost = \cost(p,  \X, \bm{S}, \bm{a}, \bm{\sigma}, \bm{K})< \infty
\end{equation}
not depending on the hidden and output layer widths $\bm{n} = (n_i)_{i=1}^L$, such that, for every $\ell \in \cur{0,1,\ldots, L}$, it holds
\begin{equation}\label{eq:induction-process} \W_p\bra{ f^{(\ell)}, G^{(\ell)}_{n_\ell}} \le \cost  n_\ell^2 \bra{\gamma_p^{|S_\ell|}}^{n_\ell} \sum_{\ell =1}^{\ell-1}\frac 1 {n_i}.\end{equation}
Moreover, for every Lipschitz continuous function $h: \R^{S_\ell} \to \R^{T}$ (with $T$ finite set) it holds
\begin{equation}
 \label{eq:induction-kernel}
\W_p\bra{  \frac 1 {n_\ell} \sum_{i=1}^{n_\ell} h(f^{(\ell)}[i,:])^{\otimes 2} , \EE\sqa{ h(G^{\ell})^{\otimes 2} } + \cN_{\Sigma_h^{(\ell)}} }   \le  \cost_h  \sum_{i=1}^{\ell} \frac{1}{n_i}
\end{equation}
where
\begin{equation}
 \cost_h =  \cost(p,  \X, \bm{S}, \bm{a}, \bm{\sigma}, \bm{K}, h, T) < \infty
\end{equation}
and the covariance $\Sigma_h^{(\ell)} \in \Sym_+^{T}$ is defined recursively by setting $\Sigma_h^{(0)} := 0$ and for $\ell\in\cur{0,1, \ldots, L}$
\begin{equation}
  \Sigma_h^{(\ell)} := \abs{D_{\cM_h}(K^{(\ell)}) \sqrt{ \tr_{a_\ell\times a_\ell} \bra{ \Sigma^{(\ell-1)}_{\sigma^{(\ell)}} } }}^2  + \frac {\cV_h( K^{(\ell)})}{n_{\ell} },
\end{equation}
with $\tr_{a_\ell \times a_\ell}$ in the sense of \eqref{eq:partial-trace-SS}.
\end{theorem}

\begin{remark}
Arguing by induction over $\ell$, using \eqref{eq:norm-cm-cv} and \eqref{eq:norm-d-cm-k}, we have the inequality
\begin{equation}\label{eq:trace-sigma-h}
 \tr( \Sigma_h^{(\ell)} ) \le \cost \nor{h}^4_{\Lip}\sum_{i=1}^{\ell} \frac{1}{n_i},
\end{equation}
where $\cost =\cost(p, \cX, \bm{k}, \bm{a}, \bm{K})<\infty$ does not depend on the layer widths $\bm{n}$. Finally, we notice the following identity between laws
\begin{equation}
 \label{eq:identity-laws-sigma-h}
 \cN_{\Sigma_h^{(\ell)}}  \stackrel{law}{=} D_{\cM_h}(K^{(\ell)}) \tr_{a_\ell} \bra{ \cN_{ \Sigma^{(\ell-1)}_{\sigma^{(\ell)}} } }+ \frac{ \cN_{\cV_h(K^{(\ell)})}}{\sqrt{n_\ell}},
\end{equation}
where the two Gaussian variable in the right hand side are independent.
\end{remark}

\begin{remark}\label{rem:large-n}
 By \Cref{thm:main-induction-basteri} and the above remark, we already know that the left hand sides of both \eqref{eq:induction-process} and \eqref{eq:induction-kernel}  are infinitesimal as the layer widths grow, hence we can (and conveniently will) restrict the proof to the case of wide enough networks.
\end{remark}

\begin{proof}
We establish the thesis arguing by induction over the layers $\ell \in \cur{0, \ldots, L}$. The case $\ell =0$ is straightforward, since $f^{(0)}$ is constant and by definition $G^{(0)} = f^{(0)}$, hence all inequalities become trivially equalities.

Next, we deduce the validity of  all the inequalities in the thesis for some $\ell \in \cur{1, \ldots, L}$ assuming that they hold for $\ell-1$ (for any $p \ge 1$ any choice of layer dimensions $\bm{n}$ and any Lipschitz function $h: \R^{S_\ell} \to \R^T$).  To simplify the notation, we write throughout the proof
 \begin{equation}
 A_0 :=  K^{(\ell)}, \quad A_1 := \frac 1 {n_{\ell-1}} \sum_{i=1}^{n_{\ell-1}} \tr_{a_\ell} \bra{ \sigma^{(\ell)} (f^{(\ell-1)} [i,:])^{\otimes 2}} \quad\text{and}\quad V_0:=  \tr_{a_\ell\times a_\ell} \bra{ \Sigma_{\sigma^{(\ell)}}^{(\ell-1)}}.
\end{equation}

First, we write \eqref{eq:induction-kernel} in the case $\ell-1$ and $h = \sigma^{(\ell)}$ which reads
\begin{equation}\label{eq:induction-square}\begin{split}
\W_p\bra{ \frac 1 {n_{\ell-1}} \sum_{i=1}^{n_{\ell-1}}  \sigma^{(\ell)} (f^{(\ell-1)} [i,:])^{\otimes 2}, \EE\sqa{ \sigma^{(\ell)}\bra{G^{(\ell-1)}}^{\otimes 2}}+  \cN_{\Sigma_{\sigma^{(\ell)}}^{(\ell-1)}}  }   &  \le  \cost_{\sigma^{(\ell)}} \sum_{i=1}^{\ell-1} \frac{1}{n_i} \\
& \le \cost \sum_{i=1}^{\ell-1} \frac{1}{n_i},
\end{split}
\end{equation}
where in the last constant $\cost$ we incorporate also the dependence upon $\sigma^{(\ell)}$. From this inequality, taking $\tr_{a_\ell}$ both sides (and using \eqref{eq:norm-partial-trace}), we obtain
 \begin{equation}\label{eq:induction-square}\begin{split}
\W_p\bra{ A_1, A_0 +  \cN_{V_0}  }    & \le \sqrt{a_{\ell}}  \W_p\bra{ \frac 1 {n_{\ell-1}} \sum_{i=1}^{n_{\ell-1}}  \sigma^{(\ell)} (f^{(\ell-1)} [i,:])^{\otimes 2}, \EE\sqa{ \sigma^{(\ell)}\bra{G^{(\ell-1)}}^{\otimes 2}}+  \cN_{\Sigma_{\sigma^{(\ell)}}^{(\ell-1)}}  }\\
& \le \cost \sum_{i=1}^{\ell-1} \frac{1}{n_i}.
\end{split}
\end{equation}
We easily deduce in particular that
\begin{equation}\label{eq:nor-a1-a0}\begin{split}
 \nor{A_1 -A_0}_{L^{p}} & = \W_{p} \bra{A_1, A_0} \\
 & \le \W_{p} \bra{ A_1, A_0 + \cN_{V_0} } + \W_p\bra{A_0 + \cN_{V_0},   A_0 }\\
 & \le \cost \sum_{i=1}^{\ell-1} \frac{1}{n_i} +  \sqrt{ \tr\bra{ V_0 }} \le \cost \sqrt{ \sum_{i=1}^{\ell-1} \frac{1}{n_i}},
 \end{split}
\end{equation}
where in the last inequality we used \eqref{eq:trace-sigma-h}.

With a little extra considerations, we further deduce that it holds
\begin{equation}\label{eq:induction-sqrt}
 \W_p\bra{  \sqrt{A_1},  \sqrt{A_0} + \cN_{\tilde{V}_0}}    \le  \cost \sum_{i=1}^{\ell-1} \frac{1}{n_i},
\end{equation}
where we set  $\tilde{V}_0 := 0$ if $\ell=1$, otherwise
\begin{equation}
\tilde{V}_0 := \abs{ D_{\sqrt{\cdot}}(K^{(\ell)}) \sqrt{ V_0 } }^2.
\end{equation}
Indeed, if $\ell=1$ we have that  \eqref{eq:induction-square} is trivially an identity, hence taking the square root both sides yields \eqref{eq:induction-sqrt}. Otherwise, we use the assumption that $A_0 = K^{(\ell)}$ is invertible,  and we  apply \Cref{prop:delta-method-sqrt}, using \eqref{eq:induction-square} and \eqref{eq:nor-a1-a0} (with $2p$ instead of $p$) to obtain \eqref{eq:induction-sqrt}. We also notice that, by \eqref{eq:norm-dsqrta}, we have
\begin{equation}\label{eq:trace-tilde-sigma}
 \tr(\tilde{V}_0) \le \cost \nor{\sqrt{ V_0 }}^2 \le \cost \sum_{i=1}^{\ell-1} \frac{1}{n_i}.
\end{equation}

Using \eqref{eq:induction-sqrt} we now prove the validity of \eqref{eq:induction-process} for $\ell$. Indeed, recalling our choice for the law of the weights $W^{(\ell)}$ in \eqref{eq:law-weights} and using \eqref{eq:identity-in-law}, we have the identity in law
\begin{equation}\label{eq:identity-law-key-proof-main}
 f^{(\ell)}  = (W^{(\ell)}\otimes \Id_{S_\ell}) \sigma^{(\ell)}(f^{(\ell-1)})\stackrel{law}{=}\bra{\Id_{n_{\ell} } \otimes \sqrt{A_1} }\cN^{(\ell)}
\end{equation}
where $\cN^{(\ell)}= \cN_{n_{\ell} \times S_{\ell} }$ is a standard Gaussian random variable with values in $\R^{n_{\ell} \times S_\ell}$ and independent of $f^{(\ell-1)}$. Similarly, by the very definition of $G^{(\ell)}_{n_\ell}$, we have
\begin{equation}
 G^{(\ell)}_{n_{\ell}} \stackrel{law}{=} \bra{\Id_{n_{\ell}} \otimes  \sqrt{ A_0  }} \cN^{(\ell)}.
\end{equation}
Hence, we collect the identity
\begin{equation}
 \W_p\bra{ f^{(\ell)}, G^{(\ell)}_{n_{\ell}} }  = \W_p\bra{ \bra{\Id_{n_{\ell} } \otimes \sqrt{A_1}}\cN^{(\ell)}, \bra{\Id_{n_{\ell}} \otimes  \sqrt{ A_0 } } \cN^{(\ell)}}.
 \end{equation}
By the triangle inequality, we split and bound from above
\begin{equation}
 \begin{split}
   \W_p\bra{ f^{(\ell)}, G^{(\ell)}_{n_{\ell}} } & \le \W_p\bra{  \bra{\Id_{n_{\ell}} \otimes \sqrt{A_1} } \cN^{(\ell)} ,  (\Id_{n_{\ell}} \otimes B) \cN^{(\ell)} }\\
  & \quad +\W_p\bra{ (\Id_{n_{\ell}} \otimes B ) \cN^{(\ell)}, \bra{\Id_{n_{\ell}} \otimes  \sqrt{ A_1 } } \cN^{(\ell)}}
  \end{split}
 \end{equation}
where we choose precisely
\begin{equation}
 B  :=  \sqrt{ A_0 } + \cN_{\tilde V_0}.
\end{equation}

By \eqref{eq:bound-abz} and \eqref{eq:induction-sqrt}, we bound the first term in the sum as follows:
  \begin{equation}\label{eq:first-term-induction}\begin{split}
    & \W_p\bra{  \bra{\Id_{n_{\ell}} \otimes  \sqrt{A_1}} \cN^{(\ell)} ,  (\Id_{n_{\ell}} \otimes B) \cN^{(\ell)} }\\
    & \le \sqrt{n_\ell} \W_p\bra{\sqrt{A_1}, B }  \le \cost \sqrt{n_\ell} \sum_{i=1}^{\ell-1} \frac{1}{n_i}.
 \end{split}
  \end{equation}
For the second term, we apply instead \Cref{prop:gaussian-closer} with
\begin{equation}
 n:= n_\ell, \quad k:= |S_\ell|, \quad B_0:= \sqrt{A_0} \quad \text{and} \quad G := \cN_{\tilde{V}_0}.
\end{equation}
Notice that, by \eqref{eq:trace-tilde-sigma} and \Cref{rem:large-n} we may assume that $\nor{ \Cov (G)} \le 1$. Then,
\begin{equation}\label{eq:second-term-induction}
 \begin{split}
  \W_p\bra{ (\Id_{n_{\ell}} \otimes B ) \cN^{(\ell)},   \bra{\Id_{n_{\ell}} \otimes \sqrt{A_0} }\cN^{(\ell)}  } & \le \cost n_\ell^2 \gamma_{p}^{|S_\ell|n_\ell } \tr\bra{ \tilde{V}_0}\\
  & \le \cost n_\ell^2 \gamma_{p}^{|S_\ell|n_\ell } \sum_{i=1}^{\ell-1} \frac 1 {n_i},
  \end{split}
 \end{equation}
 having used  \eqref{eq:trace-tilde-sigma} in the second inequality. Combining \eqref{eq:first-term-induction} with \eqref{eq:second-term-induction} yields  \eqref{eq:induction-process}.

Next, we establish \eqref{eq:induction-kernel} for $\ell$ and any Lipschitz continuous function  $h: \R^{S_{\ell}} \to \R^T$. With the current notation for $A_0$ and $A_1$, we apply \Cref{prop:clt-new}, with $n := n_\ell$ and $S := S_\ell$. Recalling \eqref{eq:identity-law-key-proof-main}, we have
\begin{equation}
 \frac 1 {n_\ell} \sum_{i=1}^{n_\ell}\bra{ h(f^{(\ell)}[i,:])}^{\otimes 2} \stackrel{law}{=} \frac 1 {n_\ell} \sum_{i=1}^{n_\ell}  h\bra{ \bra{\Id_{n_{\ell} } \otimes \sqrt{A_1} }\cN^{(\ell)} }^{\otimes 2},
\end{equation}
hence \eqref{eq:clt-quant-square} reads
\begin{equation}\begin{split}
   \W_p\bra{\frac 1 {n_\ell} \sum_{i=1}^{n_\ell}\bra{ h(f^{(\ell)}[i,:])}^{\otimes 2}, \cM_h( A_1) + \frac{\cN_{\cV_h(A_1)}}{\sqrt{n_\ell}}} & \le \cost_h \bra{ \frac 1 {n_\ell} + \nor{A_1 -A_0}_{L^{2p}}^2}\\
   & \le \cost_h \sum_{i=1}^\ell \frac 1 {n_i},
  \end{split}
\end{equation}
where  we used \eqref{eq:nor-a1-a0} (with $2p$ instead of $p$) and we stress now that $\cost$ depends on $h$ too. Therefore, in order to conclude that \eqref{eq:induction-process} holds, it is sufficient to prove that
\begin{equation}\label{eq:main-induction-process-last-step}
 \W_p\bra{\cM_h( A_1) + \frac{\cN_{\cV_h(A_1)}}{\sqrt{n_\ell}}, \cM_h(A_0) + D_{\cM_h}(A_0) \cN_{V_0} + \frac{\cN_{\cV_h(A_0)}}{\sqrt{n_\ell}}} \le \cost_h \sum_{i=1}^{\ell} \frac{1}{n_i}
\end{equation}
where $\cN_{V_0}$ and $\cN_{\cV_h(A_0)}$ are independent (recall also \eqref{eq:identity-laws-sigma-h}). We obtain this inequality in two steps. First, we have, by \eqref{eq:w-2-trivial-bound},
\begin{equation}\begin{split}\label{eq:first-step-induction-square}
 \W_p\bra{\cM_h( A_1) + \frac{\cN_{\cV_h(A_1)}}{\sqrt{n_\ell}}, \cM_h( A_1) + \frac{\cN_{\cV_h(A_0)}}{\sqrt{n_\ell}}} & \le \frac{ \nor{ \bra{ \sqrt{\cV_h(A_1)} - \sqrt{\cV_h(A_0)}}\cN_{T \times T} }_{L^p} }{\sqrt{n_\ell}}\\
 & \le \cost_T \frac{ \nor{ \sqrt{\cV_h(A_1)} - \sqrt{\cV_h(A_0)}}_{L^p}}{\sqrt{n_\ell}}\\
& \le \cost_{T,h} \frac{ \nor{ A_1 - A_0}_{L^{2p}} }{\sqrt{n_{\ell}}} \\
& \le \cost_{T,h} \frac 1 {\sqrt{n_\ell}} \sqrt{ \sum_{i=1}^{\ell-1} \frac{1}{n_i}} \le \cost_{T,h} \sum_{i=1}^{\ell} \frac{1}{n_i}.
\end{split}
\end{equation}
having used \eqref{eq:lipschitz-sqrt-cvh} in the second inequality and then \eqref{eq:nor-a1-a0} (again with $2p$ instead of $p$). Secondly, we have, by \eqref{eq:sum} and \Cref{prop:delta-method-mh},
\begin{equation}\begin{split}
 & \W_p\bra{\cM_h( A_1) + \frac{\cN_{\cV_h(A_0)}}{\sqrt{n_\ell}}, \cM_h(A_0) + D_{\cM_h}(A_0) \cN_{V_0}  + \frac{\cN_{\cV_h(A_0)}}{\sqrt{n_\ell}}} \le \\
 & \le  \W_p\bra{\cM_h( A_1), \cM_h(A_0) + D_{\cM_h}(A_0) \cN_{V_0} } \\
 & \le \cost_h \bra{ \nor{A_1-A_0}_{L^{2p}}^2 + \W_p\bra{ A_1, A_0 +  \cN_{V_0} }}\\
 & \le \cost_h \sum_{i=1}^{\ell-1} \frac 1 {n_i},
 \end{split}
\end{equation}
having used in last inequality again \eqref{eq:nor-a1-a0} and \eqref{eq:induction-square}. By the triangle inequality for $\W_p$ and using these two bounds, we conclude that \eqref{eq:main-induction-process-last-step} holds, hence \eqref{eq:induction-process} is established. \qedhere
\end{proof}

\section{Application to exact Bayesian posterior bounds}\label{sec:bayes}

In this section, we present an application of our results to bound the  exact Bayesian posteriors in a supervised learning  problem of regression type. Let us first describe a general framework to fix the notation and establish the simple abstract tool (\Cref{lem:bayes-w1}) that we then specialize in the context of deep neural networks.

In a standard regression problem, one is given a \emph{training} dataset
\begin{equation}
 \cur{ (x_i, y_i)}_{i \in \cD} \subseteq \R^{d_{in}} \times \R^{d_{out}}
\end{equation}
 and a parametrized family of functions $(f_{\theta})_{\theta \in \Theta}$, where
\begin{equation}
 f_\theta :  \R^{d_{in}} \to  \R^{d_{out}}
\end{equation}
from which one needs to estimate  $\theta^* \in \Theta$ so that $f_{\theta^*}$ best fits the dataset, according to a suitable criterion. Once $\theta^*$ is found, the resulting function $f_{\theta^*}$ can be used to compute predictions on a \emph{test} dataset $\cur{x_i}_{i \in \cT} \subseteq \R^{d_{in}}$, via $x_i \mapsto f_{\theta^*}(x_i)$. In such a generality, and without bothering too much about theoretical or computational issues, one can can roughly devise two approaches:

\begin{enumerate}
 \item a variational (or frequentist) one, where one  chooses a \emph{loss} function $\ell( f_\theta(x), y)$, e.g.\ the squared norm of the residual, $\ell( f_\theta(x), y) = \nor{ f_{\theta}(x) - y}^2$, a \emph{regularization} function $\theta\mapsto R(\theta) \in [0, \infty)$ (but one can also let $R(\theta) = 0$) and estimates $\theta^*$ by minimization of the \emph{empirical risk}:
 \begin{equation}
  \theta^* \in  \operatorname{argmin}_{ \theta } \sum_{i \in \cD} \ell(f_\theta(x_i), y_i)+ R(\theta),
 \end{equation}
 \item a probabilistic (or Bayesian) one, where one chooses a \emph{likelihood} function for the parameter $\theta$, associated to the training dataset, $\mathcal{L}(\theta; \cur{ (x_i, y_i)}_{i \in \cD}  )$, a \emph{prior} distribution of the parameter $\mathsf{p}_{prior}(\theta)$ and using Bayes' rule one obtains the (exact) \emph{posterior} distribution
 \begin{equation}
 \mathsf{p}_{posterior} (\theta | \cur{ (x_i, y_i)}_{i \in \cD}) \propto \mathcal{L}(\theta; \cur{ (x_i, y_i)}_{i \in \cD}  ) \mathsf{p}_{prior}(\theta),
 \end{equation}
 where $\propto$ means that the identity holds true up to a multiplicative constant ensuring that the left-hand side is indeed a probability distribution. A possible criterion to define $\theta^*$ is then the \emph{maximum a posteriori} estimate:
 \begin{equation}\begin{split}
  \theta^* &\in  \operatorname{argmax}_{ \theta } \mathsf{p}_{posterior} (\theta | \cur{ (x_i, y_i)}_{i \in \cD}) \\
  & \quad  = \operatorname{argmax}_{ \theta }\mathcal{L}(\theta; \cur{ (x_i, y_i)}_{i \in \cD}  ) \mathsf{p}_{prior}(\theta).
  \end{split}
 \end{equation}
\end{enumerate}

The two approaches are evidently related, and a correspondence can be obtained by taking a logarithm and changing  sign in the Bayesian one, or equivalently, starting from the variational one, for a given a loss function $\ell$, one defines as a likelihood the function
\begin{equation}
 \mathcal{L}(\theta; \cur{ (x_i, y_i)}_{i \in \cD}  ) := \exp\bra{ - \sum_{i\in \cD} \ell( f_\theta(x_i), y_i)}
\end{equation}
and for a given regularization $R$, one defines the prior distribution
\begin{equation}
 \mathsf{p}_{prior}(\theta) \propto \exp\bra{ - R(\theta)}.
\end{equation}
With this construction, it is immediate to see that the estimated $\theta^*$'s  are the same in the two approaches. Let us notice that in the example of the squared norm of the residual we obtain a Gaussian likelihood
\begin{equation}
 \mathcal{L}(\theta; \cur{ (x_i, y_i)}_{i \in \cD}  )  = \exp\bra{ - \sum_{i \in \cD} \nor{ f_\theta(x_i) - y_i}^2 },
\end{equation}
which has the probabilistic interpretation that the residuals $f_{\theta}(x_i) -y_i$ are independent Gaussian variables, conditionally upon the knowledge of $\theta \in \Theta$.

The general question that we consider can be phrased as follows: given two different prior distributions, that are close in a suitable metric, can one provide bounds on the associated posteriors (and possibly on related quantities, such as $\theta^*$)? The answer clearly depends on the chosen metric, and in our case it is natural to consider the Wasserstein distance, which however seems too weak to provide meaningful bounds on the estimated $\theta^*$.

Precisely, we make use of the following general result.

\begin{lemma}\label{lem:bayes-w1}
 Let $\mu$, $\nu$ be probability measures on $\R^S$ for some finite set $S$ and finite absolute moments of order $p \ge 1$. Write for $q \ge 1$,
 \begin{equation}
  m_{q}(\mu) := \int_{\R^S} \nor{z}^{q} d \mu(z).
 \end{equation}
 Let $g: \R^S \to [0, \nor{g}_0]$ be Lipschitz continuous and uniformly bounded (i.e., of polynomial growth of order $0$), such that
 \begin{equation}
  \mu(g):= \int g d \mu >0 \quad \text{and} \quad \nu(g):= \int g d \nu >0,
 \end{equation}
 and define the probability measures
 \begin{equation}
  \mu_g :=\frac{ g }{\mu(g)} \mu \quad \nu_g := \frac{ g }{\nu(g) } \nu.
 \end{equation}
 Then,
 \begin{equation}
  \W_1(\mu_g, \nu_g) \le \frac {1}{\mu(g)} \bra{  \Lip(g) m_{p'}(\mu) +   \bra{ 1 + \frac{ m_1(\mu) \Lip(g)}{\nu(g)} }\nor{g}_0  }\W_p(\mu, \nu),
 \end{equation}
  where $p' = p/(p-1)$.
 \end{lemma}

 \begin{proof}
 Let us first notice that, given (jointly defined) random variables $X$ and $Y$ with laws $\mathbb{P}_X = \mu$, $\mathbb{P}_Y = \nu$, we have
 \begin{equation}
  \abs{ \mu(g) - \nu(g)} = \abs{ \EE\sqa{ g(X) - g(Y)} } \le \nor{ g(X) - g(Y)}_{L^1} \le \Lip(g) \nor{X-Y}_{L^1},
 \end{equation}
 hence minimizing upon the joint realizations of $X$ and $Y$ we have
 \begin{equation}\label{eq:bound-mug-nug}
  \abs{ \mu(g) - \nu(g)} \le \Lip(g) \W_1(\mu, \nu).
 \end{equation}
 Similarly, given a function $f: \R^S \to \R$ with $\Lip(f) \le 1$ and $f(0) =0$, so that $\abs{f(z)} \le \nor{z}$ for every $z \in \R^S$, we find
 \begin{equation}\begin{split}
  \abs{ \int f g d \mu -  \int f g\nu } & \le \nor{f(X)g(X) - f(Y)g(Y)}_{L^1}  \\
  & \le \nor{ f(X) (g(X)-g(Y))}_{L^1} +  \nor{ (f(X)-f(Y))g(Y)}_{L^1} \\
  & \le \nor{f(X)}_{L^{p'}}  \nor{ g(X)-g(Y)}_{L^p} + \nor{ f(X)-f(Y)}_{L^1} \nor{g(Y)}_{L^\infty}\\
  & \le m_{p'}(\mu) \Lip(g) \nor{ X-Y}_{L^p}  + \nor{g}_0 \nor{ X-Y}_{L^p},
  \end{split}
 \end{equation}
hence minimization  upon the joint realizations of $X$ and $Y$ gives
 \begin{equation}\label{eq:bound-mufg-nufg}
\abs{ \int f g d \mu -  \int f g\nu } \le \bra{ m_{p'}(\mu) \Lip(g) + \nor{g}_0 } \W_p\bra{X,Y}.
 \end{equation}

 Next, we use  the standard Kantorovich duality for the $\W_1$ distance (see \cite{villani2009optimal}), that yields
\begin{equation}
  \W_1(\mu_g, \nu_g) = \sup_{\substack{f: \R^S \to \R\\  \Lip(f) \le 1 }} \cur{ \int f d \mu_g - \int f d \nu_g}.
\end{equation}
Since $f$ in the maximization above can be replaced by $f + \cost$ for any constant $\cost< \infty$, we can assume without loss of generality that $f(0) = 0$, hence $|f(z) | \le \nor{z}$ for every $z \in \R^S$. Then,
\begin{equation} \begin{split} \int f d \mu_g - \int f d \nu_g &  = \frac{1}{\mu(g)}  \int f g d (\mu - \nu) + \bra{\frac{1}{\mu(g)} - \frac{1}{\nu(g)}} \int fg d \nu\\
  & \le  \frac{1}{\mu(g)}\bra{ \int f g d \mu - \int fg \nu} +   \nor{g}_0 \frac{ \abs{\mu(g) - \nu(g)}}{ \mu(g)\nu(g)} \int \|x\|  d \nu\\
  & \le \frac{\bra{ m_{p'}(\mu) \Lip(g) + \nor{g}_0 } }{\mu(g)} \W_p\bra{X,Y} + \frac{ \Lip(g)  \nor{g}_0 m_1(\mu)}{\mu(g)\nu(g)}\W_1(\mu, \nu), \end{split}
\end{equation}
having used \eqref{eq:bound-mufg-nufg}  and \eqref{eq:bound-mug-nug}. Since $\W_1(\mu, \nu) \le \W_p(\mu, \nu)$, the thesis follows.
 \end{proof}

 \begin{remark}
  One can extend the previous argument to other distances defined by maximization over a suitable family of ``test'' functions: for example,  in the case of the bounded-Lipschitz distance (also called Dudley metric)
  \begin{equation}
   d_{BL}(\mu, \nu) := \sup_{ \substack{ \Lip(f)\le 1 \\ \nor{f}_0 \le 1 } } \cur{ \int f d \mu - \int f d \nu},
  \end{equation}
  we may bound from above $d_{BL}(\mu_g, \nu_g)$ in terms of $\W_1(\mu, \nu)$.
 \end{remark}

We next apply the general result in the context of deep neural networks with random Gaussian weights and their Gaussian  approximations. Deep neural networks $f^{(L)}$ appear to be naturally a parametrized family of functions with $\theta := \bm{W}$ (using the generalized architecture to treat weights and biases at once) with $d_{out} = n_L$. On the other side, their Gaussian approximation $G^{(L)}_{n_L}$ may be thought as a parametrized family letting  $\theta := G^{(L)}$ itself. \Cref{thm:main-induction-basteri} and \Cref{thm:main-induction} can be interpreted as measuring the distance between two \emph{prior} distributions. Hence, by \Cref{lem:bayes-w1} we are able to obtain a result on the associated exact \emph{posterior} distributions given any training dataset and a test set (both finite).

The key assumption that we impose is that the likelihood function can be written in the form
\begin{equation}\label{eq:likelihood-g}
\mathcal{L} ( \theta;  \cur{(x_i,y_i)}_{i\in \cD}) :=  g\bra{  (f_\theta(x_i))_{i \in \cD} }
\end{equation}
for some $g: \R^{d_{out} \times \cD} \to [0, \infty)$ that is Lipschitz and uniformly bounded (the exact dependence upon $(y_i)_{i\in \cD}$ is not relevant). We can then use $\mathcal{L}$ to define a notion of posterior distribution: if $F_{train}$ is any random variable defined on some probability space $(\Omega, \mathcal{A}, \mathbb{P})$ with values in $\mathbb{R}^{d_{out} \times \cD}$, then the posterior distribution given the dataset is abstractly given by
\begin{equation}
\mathbb{P}_{|\cD} :=  \frac{ g(F_{train})}{\int g(F_{train}) d\mathbb{P}} \mathbb{P},
\end{equation}
provided that the denominator is strictly positive. Let us also notice that, given the specific form of the likelihood \eqref{eq:likelihood-g}, by an abstract change of variables one can reduce to any ``intermediate'' space between $\Omega$ and $\R^{d_{out} \times \cD}$. In particular, if one is only interested in considering a further random variable $F_{test}: \Omega \to \R^{d_{out} \times \cT}$, we can  reduce the considerations from an ``abstract'' $\Omega$ to $\R^{d_{out} \times (\cD \cup \cT)}$, endowed with the ``prior'' distribution given by the joint law of $F := (F_{train}, F_{test})$ with respect to $\mathbb{P}$. But then, to keep the notation simple, we can further assume that $g$ is defined on the entire domain $\R^{d_{out} \times (\cD \cup \cT)}$ and denote with $F_{|\cD}$ the law of $F$ with respect to $\mathbb{P}_{|\cD}$. In view of these considerations, it is rather straightforward to obtain the following result.

\begin{corollary}\label{cor:posterior}
Consider a training dataset $\cur{(x_i, y_i)}_{i \in \cD}$ and a test set $\cur{x_j}_{j \in \cT}$ (both finite) and define $\cX := \cur{x_i}_{i\in \cD} \cup \cur{x_j}_{j \in \cT}$.  With the notation of \Cref{sec:nn},  assume that all the activation functions of a (generalized) deep neural network architecture $\bm{\sigma} = (\sigma^{(\ell)})_{\ell = 1}^L$ are Lipschitz continuous and  the weights $\bm{W}$ are random Gaussian variables as in \eqref{eq:law-weights}.

Consider any likelihood function of the form \eqref{eq:likelihood-g} with $g: \R^{d_{out} \times S_L} \to [0, \infty)$ that is Lipschitz continuous, uniformly bounded and such that $\EE\sqa{ g( G^{(L)}_{n_L} )} >0$.

Then, there exists
\begin{equation}
  \cost = \cost(\X, \bm{S}, \bm{a}, \bm{\sigma}, \bm{K}, g, n_L)< \infty
\end{equation}
such that, if $n: = \min_{i=1, \ldots, L-1} \cur{n_i} \ge \cost$, then $\EE\sqa{ g( f^{(L)} )}>0$, the posterior laws of $f^{(L)}$ and $G^{(L)}_{n_L}$ are both well defined and it holds
\begin{equation} \W_1\bra{   f^{(L)}_{|\cD} , G^{(L)}_{n_L, | \cD}  }  \le \cost \frac 1 {\sqrt{n}}. \end{equation}
If moreover the infinite-width covariance matrices $\bm{K}$ are non-degenerate on $\cX$, one also has
\begin{equation}\label{eq:posterior-non-degenerate} \W_1\bra{   f^{(L)}_{|\cD} , G^{(L)}_{n_L | \cD}  }  \le \cost \frac 1 n. \end{equation}
(possibly with a different constant $\cost < \infty$).
\end{corollary}

\begin{remark}
 Notice that in order for  \eqref{eq:posterior-non-degenerate} to hold, the non-degeneracy assumption must hold on the entire input set $\cX$, which includes both the training and the test (input) sets.
\end{remark}

\begin{proof}

 It is sufficient to apply \Cref{lem:bayes-w1} with $p=2$ and $\mu$, $\nu$ given respectively by the laws of  $G^{(L)}_{n_L}$ and $f^{(L)}$  on the input set $\cX$. Indeed, we notice that by \eqref{eq:bound-mug-nug} it holds
 \begin{equation}
  \abs{ \EE\sqa{ g  \bra{ G^{(L)}_{n_L}}} - \EE\sqa{ g\bra{ f^{(L)} }} } \le \cost \W_1 \bra{  f^{(L)},  G^{(L)}_{n_L} },
 \end{equation}
 hence if $n$ is large, using \Cref{thm:main-induction-basteri} we see that also the term $\EE\sqa{ g( f^{(L)} )}$ will be positive, and actually larger than $\EE\sqa{ g  \bra{ G^{(L)}_{n_L}}} /2$. Using this fact we only need to rephrase \Cref{lem:bayes-w1} using the probabilistic notation and notice that the absolute moments $m_1(\mu)$ and $m_2(\mu)$ are obviously finite because $G^{(L)}_{n_L}$ is Gaussian.
\end{proof}

\begin{remark}
In the case of the squared residual loss function, we have already commented that that the corresponding likelihood is Gaussian, i.e.,
\begin{equation}
 \mathcal{L}(\theta; \cD) = \exp\bra{ - \sum_{i \in \cD} \| f_\theta(x_i) - y_i \|^2 } = g\bra{ (f_\theta(x_i))_{i \in \cD} },
\end{equation}
where
 \begin{equation}
  g( (z_i)_{i \in \cD}) := \exp\bra{ - \sum_{i \in \cD}  \| z_i - y_i \|^2}.
 \end{equation}
Clearly $g$ is Lipschitz continuous (with a Lipschitz constant $\Lip(g)$ that depends on the size of the training set $\cD$), uniformly bounded (by the constant $1$) and positive everywhere. It is also simple to argue that the posterior distribution of $G^{(L)}_{n_L | \cD}$ is a Gaussian law (the mean and covariance parameters can be also written, see e.g.\ \cite{williams2006gaussian} for more on the subject). Hence, we  deduce that also the exact Bayesian posterior distribution of a deep neural network (evaluated on a finite test set and with the Gaussian prior distribution on the weights \eqref{eq:law-weights}) can be quantitatively approximated with a Gaussian random variable with a rate $1/\sqrt{n}$ (or $1/n$) where  $n: =\min_{i=1, \ldots, L-1}\cur{n_i}$, depending on the (non-)degeneracy of the infinite-width covariances.
\end{remark}

\printbibliography

\end{document}